%
\magnification1200
\baselineskip15pt
\hfuzz2pt
\newcount\relFnno
\def\ref#1{\expandafter\edef\csname#1\endcsname}
\relFnno 1
\ref {Introduction}{1}
\ref {E50}{$(1.2)$}
\relFnno 1
\ref {E61}{$(1.6)$}
\ref {Fn1}{\the \relFnno }
\advance \relFnno 1
\relFnno 1
\ref {E60}{$(1.11)$}
\ref {E62}{$(1.12)$}
\ref {Fn2}{\the \relFnno }
\advance \relFnno 1
\relFnno 1
\relFnno 1
\ref {Fourier}{2}
\relFnno 1
\ref {E2}{$(2.5)$}
\ref {HarishChandra}{Theorem\penalty 10000\ 2.1}
\ref {Fourier-transform}{Theorem\penalty 10000\ 2.2}
\ref {E1}{$(2.6)$}
\ref {E3}{$(2.7)$}
\relFnno 1
\ref {E5}{$(2.12)$}
\ref {E4}{$(2.13)$}
\ref {E6}{$(2.14)$}
\relFnno 1
\ref {E7}{$(2.15)$}
\ref {Capelli}{3}
\ref {Fn3}{\the \relFnno }
\advance \relFnno 1
\relFnno 1
\relFnno 1
\ref {Main-Prop-E}{Theorem\penalty 10000\ 3.1}
\ref {ExtraVanishing}{Theorem\penalty 10000\ 3.2}
\ref {Def-f-l}{$(3.6)$}
\ref {Cut-Off}{Lemma\penalty 10000\ 3.3}
\relFnno 1
\ref {nilp}{Theorem\penalty 10000\ 3.4}
\ref {roughstructure}{Lemma\penalty 10000\ 3.6}
\relFnno 1
\ref {E30}{$(3.13)$}
\ref {E51}{$(3.14)$}
\ref {Pieri1}{Theorem\penalty 10000\ 3.7}
\ref {P-Formula-1}{$(3.16)$}
\ref {Fn4}{\the \relFnno }
\advance \relFnno 1
\relFnno 1
\ref {Pieriell}{$(3.17)$}
\ref {Duality}{4}
\ref {FiltIso}{Proposition\penalty 10000\ 4.1}
\relFnno 1
\ref {AAA}{Corollary\penalty 10000\ 4.2}
\ref {Dual2}{Theorem\penalty 10000\ 4.3}
\ref {DualFormula}{$(4.6)$}
\ref {E9}{$(4.7)$}
\ref {Fn5}{\the \relFnno }
\advance \relFnno 1
\ref {Fn6}{\the \relFnno }
\advance \relFnno 1
\relFnno 1
\ref {Symm}{Corollary\penalty 10000\ 4.4}
\ref {Symmetrisch}{$(4.11)$}
\ref {Pieri2}{Theorem\penalty 10000\ 4.5}
\ref {pieriII}{$(4.12)$}
\ref {Fn7}{\the \relFnno }
\advance \relFnno 1
\relFnno 1
\ref {EvalFormula}{Corollary\penalty 10000\ 4.6}
\ref {E8}{$(4.15)$}
\ref {E52}{$(4.16)$}
\ref {BBB}{Corollary\penalty 10000\ 4.7}
\ref {virdim}{Theorem\penalty 10000\ 4.8}
\relFnno 1
\ref {Interpolation}{5}
\ref {Involution}{Theorem\penalty 10000\ 5.1}
\ref {E41}{$(5.1)$}
\ref {Interpol}{Theorem\penalty 10000\ 5.2}
\ref {E11}{$(5.3)$}
\relFnno 1
\ref {E13}{$(5.5)$}
\ref {E10}{$(5.6)$}
\ref {InvAut}{Theorem\penalty 10000\ 5.3}
\ref {E12}{$(5.8)$}
\ref {Fn8}{\the \relFnno }
\advance \relFnno 1
\relFnno 1
\ref {Product}{6}
\ref {E16}{$(6.1)$}
\ref {E15}{$(6.2)$}
\ref {E42}{$(6.5)$}
\relFnno 1
\ref {E14}{$(6.6)$}
\ref {Adjoint1}{Theorem\penalty 10000\ 6.3}
\ref {sl-two}{7}
\relFnno 1
\ref {adjoints}{Theorem\penalty 10000\ 7.1}
\ref {E17}{$(7.2)$}
\ref {triple}{Theorem\penalty 10000\ 7.2}
\relFnno 1
\ref {integration}{Theorem\penalty 10000\ 7.3}
\ref {effect}{Proposition\penalty 10000\ 7.4}
\ref {LL*}{Theorem\penalty 10000\ 7.5}
\relFnno 1
\ref {L-Formula}{$(7.11)$}
\ref {L*-Formula}{$(7.12)$}
\ref {L*-Formula2}{$(7.13)$}
\ref {E20}{$(7.14)$}
\ref {E21}{$(7.15)$}
\ref {E22}{$(7.17)$}
\relFnno 1
\ref {minusaction}{Theorem\penalty 10000\ 7.6}
\relFnno 1
\ref {differential}{8}
\ref {diffop}{Proposition\penalty 10000\ 8.1}
\relFnno 1
\relFnno 1
\ref {E23}{$(8.9)$}
\ref {difflim}{Theorem\penalty 10000\ 8.3}
\relFnno 1
\ref {E25}{$(8.14)$}
\ref {E24}{$(8.15)$}
\relFnno 1
\ref {Fn9}{\the \relFnno }
\advance \relFnno 1
\relFnno 1
\ref {Binomial}{9}
\ref {BinomialTheorem}{Theorem\penalty 10000\ 9.1}
\ref {BinomialFormula}{$(9.5)$}
\ref {Fn10}{\the \relFnno }
\advance \relFnno 1
\ref {Fn11}{\the \relFnno }
\advance \relFnno 1
\relFnno 1
\ref {DualFormula2}{$(9.9)$}
\ref {DualFormula3}{$(9.10)$}
\ref {order}{$(9.14)$}
\relFnno 1
\ref {Beispiel}{10}
\relFnno 1
\ref {E40}{$(10.8)$}
\ref {E43}{$(10.10)$}
\relFnno 1
\ref {E44}{$(10.12)$}
\ref {E45}{$(10.13)$}
\relFnno 1
\relFnno 1
\ref {References}{11}
\ref {BenRat1}{[BR1]}
\ref {BenRat}{[BR2]}
\ref {Deb}{[De]}
\ref {HoUm}{[HU]}
\ref {Las}{[La1]}
\ref {La2}{[La2]}
\ref {Kac}{[Kac]}
\ref {Annals}{[Kn1]}
\ref {Montreal}{[Kn2]}
\ref {Semi}{[Kn3]}
\ref {Existence}{[Kn4]}
\ref {SymCap}{[KnSa]}
\ref {Leahy}{[Le]}
\ref {Ok}{[Ok]}
\ref {OO1}{[OO1]}
\ref {OlOk}{[OO2]}
\ref {Sahi}{[Sa]}
\ref {Sek}{[Se]}
\ref {Up}{[Up]}
\ref {Yan}{[Yan]}
\ref {Spaltenbreite}{32.77786pt}
\relFnno 1

\font\sevenex=cmex7
\scriptfont3=\sevenex
\font\fiveex=cmex10 scaled 500
\scriptscriptfont3=\fiveex

\def\XS{{\widetilde X}}

\def\phi{\varphi}
\def\epsilon{\varepsilon}
\def\theta{\vartheta}
\def\uauf{\lower1.7pt\hbox to 3pt{%
\vbox{\offinterlineskip
\hbox{\vbox to 8.5pt{\leaders\vrule width0.2pt\vfill}%
\kern-.3pt\hbox{\lams\char"76}\kern-0.3pt%
$\raise1pt\hbox{\lams\char"76}$}}\hfil}}

\def\title#1{\par
{\baselineskip1.5\baselineskip\rightskip0pt plus 5truecm
\leavevmode\vskip0truecm\noindent\font\BF=cmbx10 scaled \magstep2\BF #1\par}
\vskip1truecm
\leftline{\font\CSC=cmcsc10\CSC Friedrich Knop}
\leftline{Department of Mathematics, Rutgers University, New Brunswick NJ
08903, USA}
\leftline{knop@math.rutgers.edu}
\vskip1truecm
\par}

\def\cite#1{\expandafter\ifx\csname#1\endcsname\relax
{\bf?}\immediate\write16{#1 ist nicht definiert!}\else\csname#1\endcsname\fi}
\def\expandwrite#1#2{\edef\next{\write#1{#2}}\next}
\def\neverexpand{\noexpand\noexpand\noexpand}
\def\strip#1\ {}
\def\ncite#1{\expandafter\ifx\csname#1\endcsname\relax
{\bf?}\immediate\write16{#1 ist nicht definiert!}\else
\expandafter\expandafter\expandafter\strip\csname#1\endcsname\fi}
\newwrite\AUX
\immediate\openout\AUX=\jobname.aux
\font\eightrm=cmr8\font\sixrm=cmr6
\font\eighti=cmmi8
\font\eightit=cmti8
\font\eightbf=cmbx8
\font\eightcsc=cmcsc10 scaled 833
\def\eightpoint{%
\textfont0=\eightrm\scriptfont0=\sixrm\def\rm{\fam0\eightrm}%
\textfont1=\eighti
\textfont\bffam=\eightbf\def\bf{\fam\bffam\eightbf}%
\textfont\itfam=\eightit\def\it{\fam\itfam\eightit}%
\def\csc{\eightcsc}%
\setbox\strutbox=\hbox{\vrule height7pt depth2pt width0pt}%
\normalbaselineskip=0,8\normalbaselineskip\normalbaselines\rm}
\newcount\absFnno\absFnno1
\write\AUX{\relFnno1}
\newif\ifMARKE\MARKEtrue
{\catcode`\@=11
\gdef\footnote{\ifMARKE\edef\@sf{\spacefactor\the\spacefactor}\/%
$^{\cite{Fn\the\absFnno}}$\@sf\fi
\MARKEtrue
\insert\footins\bgroup\eightpoint
\interlinepenalty100\let\par=\endgraf
\leftskip=0pt\rightskip=0pt
\splittopskip=10pt plus 1pt minus 1pt \floatingpenalty=20000\smallskip
\item{$^{\cite{Fn\the\absFnno}}$}%
\expandwrite\AUX{\neverexpand\ref{Fn\the\absFnno}{\neverexpand\the\relFnno}}%
\global\advance\absFnno1\write\AUX{\advance\relFnno1}%
\bgroup\strut\aftergroup\@foot\let\next}}
\skip\footins=12pt plus 2pt minus 4pt
\dimen\footins=30pc
\output={\plainoutput\immediate\write\AUX{\relFnno1}}
\newcount\Abschnitt\Abschnitt0
\def\beginsection#1. #2 \par{\advance\Abschnitt1%
\vskip0pt plus.10\vsize\penalty-250
\vskip0pt plus-.10\vsize\bigskip\vskip\parskip
\edef\TEST{\number\Abschnitt}
\expandafter\ifx\csname#1\endcsname\TEST\relax\else
\immediate\write16{#1 hat sich geaendert!}\fi
\expandwrite\AUX{\neverexpand\ref{#1}{\TEST}}
\leftline{\marginnote{#1}\bf\number\Abschnitt. \ignorespaces#2}%
\nobreak\smallskip\noindent\SATZ1\GNo0}
\def\Proof:{\par\noindent{\it Proof:}}
\def\Remark:{\ifdim\lastskip<\medskipamount\removelastskip\medskip\fi
\noindent{\bf Remark:}}
\def\Remarks:{\ifdim\lastskip<\medskipamount\removelastskip\medskip\fi
\noindent{\bf Remarks:}}
\def\Definition:{\ifdim\lastskip<\medskipamount\removelastskip\medskip\fi
\noindent{\bf Definition:}}
\def\Example:{\ifdim\lastskip<\medskipamount\removelastskip\medskip\fi
\noindent{\bf Example:}}
\def\Examples:{\ifdim\lastskip<\medskipamount\removelastskip\medskip\fi
\noindent{\bf Examples:}}
\newif\ifmarginalnotes

\def\marginnote#1{\ifmarginalnotes\hbox to 0pt{\eightpoint\hss #1\ }\fi}

\newcount\SATZ\SATZ1
\def\proclaim #1. #2\par{\ifdim\lastskip<\medskipamount\removelastskip
\medskip\fi
\noindent{\bf#1.\ }{\it#2}\Par
\ifdim\lastskip<\medskipamount\removelastskip\goodbreak\medskip\fi}
\def\Aussage#1{%
\expandafter\def\csname#1\endcsname##1.{\ifx?##1?\relax\else
\edef\TEST{#1\penalty10000\ \number\Abschnitt.\number\SATZ}
\expandafter\ifx\csname##1\endcsname\TEST\relax\else
\immediate\write16{##1 hat sich geaendert!}\fi
\expandwrite\AUX{\neverexpand\ref{##1}{\TEST}}\fi
\proclaim {\marginnote{##1}\number\Abschnitt.\number\SATZ. #1\global\advance\SATZ1}.}}
\Aussage{Theorem}
\Aussage{Proposition}
\Aussage{Corollary}
\Aussage{Lemma}
\font\la=lasy10
\def\strich{\hbox{$\vcenter{\hbox
to 1pt{\leaders\hrule height -0,2pt depth 0,6pt\hfil}}$}}
\def\dashedrightarrow{\hbox{%
\hbox to 0,5cm{\leaders\hbox to 2pt{\hfil\strich\hfil}\hfil}%
\kern-2pt\hbox{\la\char\string"29}}}

\def\Bindestrich{\penalty10000-\hskip0pt}
\let\_=\Bindestrich
\def\.{{\sfcode`.=1000.}}

\def\Par{\par}
\def\:={\mathrel{\raise0,9pt\hbox{.}\kern-2,77779pt
\raise3pt\hbox{.}\kern-2,5pt=}}
\def\=:{\mathrel{=\kern-2,5pt\raise0,9pt\hbox{.}\kern-2,77779pt
\raise3pt\hbox{.}}} 
\def\into{\hookrightarrow}
\def\pfeil{\rightarrow}

\def\pf#1{\buildrel#1\over\rightarrow}
\def\Pf#1{\buildrel#1\over\longrightarrow}

\def\Ugleich{\hbox{$\cup$\kern.5pt\vrule depth -0.5pt}}
\def\|#1|{\mathop{\rm#1}\nolimits}
\def\<{\langle}
\def\>{\rangle}
\let\Times=\times
\def\times{\mathop{\Times}}
\let\Otimes=\otimes
\def\otimes{\mathop{\Otimes}}
\catcode`\@=11
\def\hex#1{\ifcase#1 0\or1\or2\or3\or4\or5\or6\or7\or8\or9\or A\or B\or
C\or D\or E\or F\else\message{Warnung: Setze hex#1=0}0\fi}
\def\fontdef#1:#2,#3,#4.{%
\alloc@8\fam\chardef\sixt@@n\FAM
\ifx!#2!\else\expandafter\font\csname text#1\endcsname=#2
\textfont\the\FAM=\csname text#1\endcsname\fi
\ifx!#3!\else\expandafter\font\csname script#1\endcsname=#3
\scriptfont\the\FAM=\csname script#1\endcsname\fi
\ifx!#4!\else\expandafter\font\csname scriptscript#1\endcsname=#4
\scriptscriptfont\the\FAM=\csname scriptscript#1\endcsname\fi
\expandafter\edef\csname #1\endcsname{\fam\the\FAM\csname text#1\endcsname}
\expandafter\edef\csname hex#1fam\endcsname{\hex\FAM}}
\catcode`\@=12 

\fontdef Ss:cmss10,,.
\fontdef Fr:eufm10,eufm7,eufm5.

\def\fg{{\Fr g}}


\def\fn{{\Fr n}}

\def\fp{{\Fr p}}

\def\fs{{\Fr s}}
\def\ft{{\Fr t}}
\def\fU{{\Fr U}}

\def\fZ{{\Fr Z}}
\fontdef bbb:msbm10,msbm7,msbm5.
\fontdef mbf:cmmib10,cmmib7,.
\fontdef msa:msam10,msam7,msam5.
\def\CC{{\bbb C}}

\def\NN{{\bbb N}}\def\OO{{\bbb O}}

\def\ZZ{{\bbb Z}}
\def\cA{{\cal A}}\def\cB{{\cal B}}\def\cC{{\cal C}}\def\cD{{\cal D}}

\def\cK{{\cal K}}
\def\cP{{\cal P}}

\def\cZ{{\cal Z}}
\mathchardef\leer=\string"0\hexbbbfam3F
\mathchardef\subsetneq=\string"3\hexbbbfam24
\mathchardef\semidir=\string"2\hexbbbfam6E
\mathchardef\dirsemi=\string"2\hexbbbfam6F
\mathchardef\haken=\string"2\hexmsafam78
\mathchardef\auf=\string"3\hexmsafam10
\let\OL=\overline
\def\overline#1{{\hskip1pt\OL{\hskip-1pt#1\hskip-.3pt}\hskip.3pt}}

\def\cq{{\overline{c}}}
\def\Dq{{\overline{D}}}
\def\Eq{{\overline{E}}}

\def\hq{{\overline{h}}}

\def\Lq{{\overline{L}}}
\def\Mq{{\overline{M}}}

\def\pq{{\overline{p}}}
\def\qq{{\overline{q}}}
\def\Rq{{\overline{R}}}

\def\Wq{{\overline{W}}}
\def\Xq{{\overline{X}}}

%
\abovedisplayskip 9.0pt plus 3.0pt minus 3.0pt
\belowdisplayskip 9.0pt plus 3.0pt minus 3.0pt
\newdimen\Grenze\Grenze2\parindent\advance\Grenze1em
\newdimen\Breite
\newbox\DpBox
\def\NewDisplay#1$${\Breite\hsize\advance\Breite-\hangindent
\setbox\DpBox=\hbox{\hskip2\parindent$\displaystyle{#1}$}%
\ifnum\predisplaysize<\Grenze\abovedisplayskip\abovedisplayshortskip
\belowdisplayskip\belowdisplayshortskip\fi
\global\futurelet\nexttok\WEITER}
\def\WEITER{\ifx\nexttok\qed\expandafter\leftQEDdisplay
\else\leftdisplay\fi}
\def\leftdisplay{\hskip-\hangindent\leftline{\box\DpBox}$$}
\def\leftQEDdisplay{\hskip-\hangindent
\line{\copy\DpBox\hfill\lower\dp\DpBox\copy\QEDbox}%
\belowdisplayskip0pt$$\bigskip\let\nexttok=}
\everydisplay{\NewDisplay}
\newcount\GNo\GNo=0
\def\eqno#1{
\global\advance\GNo1
\edef\FTEST{$(\number\Abschnitt.\number\GNo)$}
\ifx?#1?\relax\else
\expandafter\ifx\csname#1\endcsname\FTEST\relax\else
\immediate\write16{#1 hat sich geaendert!}\fi
\expandwrite\AUX{\neverexpand\ref{#1}{\FTEST}}\fi
\llap{\hbox to 40pt{\marginnote{#1}\FTEST\hfill}}}

\catcode`@=11
\def\eqalignno#1{\null\vcenter{\openup\jot\m@th\ialign{\eqno{##}\hfil
&\strut\hfil$\displaystyle{##}$&$\displaystyle{{}##}$\hfil\crcr#1\crcr}}\,}
\catcode`@=12

\newbox\QEDbox
\newbox\nichts\setbox\nichts=\vbox{}\wd\nichts=2mm\ht\nichts=2mm
\setbox\QEDbox=\hbox{\vrule\vbox{\hrule\copy\nichts\hrule}\vrule}
\def\qed{\leavevmode\unskip\hfil\null\nobreak\hfill\copy\QEDbox\medbreak}
\newdimen\HIindent
\newbox\HIbox
\def\setHI#1{\setbox\HIbox=\hbox{#1}\HIindent=\wd\HIbox}
\def\HI#1{\par\hangindent\HIindent\hangafter=0\noindent\leavevmode
\llap{\hbox to\HIindent{#1\hfil}}\ignorespaces}

\newdimen\maxSpalbr
\newdimen\altSpalbr
\newcount\Zaehler

\def\beginrefs{%
\expandafter\ifx\csname Spaltenbreite\endcsname\relax
\def\Spaltenbreite{1cm}\immediate\write16{Spaltenbreite undefiniert!}\fi
\expandafter\altSpalbr\Spaltenbreite
\maxSpalbr0pt
\gdef\alt{}
\def\\##1\relax{%
\gdef\neu{##1}\ifx\alt\neu\global\advance\Zaehler1\else
\xdef\alt{\neu}\global\Zaehler=1\fi\xdef\SigText{##1\the\Zaehler}}
\def\L|Abk:##1|Sig:##2|Au:##3|Tit:##4|Zs:##5|Bd:##6|S:##7|J:##8||{%
\def\SigText{##2}\global\setbox0=\hbox{##2\relax}
\edef\TEST{[\SigText]}
\expandafter\ifx\csname##1\endcsname\TEST\relax\else
\immediate\write16{##1 hat sich geaendert!}\fi
\expandwrite\AUX{\neverexpand\ref{##1}{\TEST}}
\setHI{[\SigText]\ }
\ifnum\HIindent>\maxSpalbr\maxSpalbr\HIindent\fi
\ifnum\HIindent<\altSpalbr\HIindent\altSpalbr\fi
\HI{\marginnote{##1}[\SigText]}
\ifx-##3\relax\else{##3}: \fi
\ifx-##4\relax\else{##4}{\sfcode`.=3000.} \fi
\ifx-##5\relax\else{\it ##5\/} \fi
\ifx-##6\relax\else{\bf ##6} \fi
\ifx-##8\relax\else({##8})\fi
\ifx-##7\relax\else, {##7}\fi\Par}
\def\B|Abk:##1|Sig:##2|Au:##3|Tit:##4|Reihe:##5|Verlag:##6|Ort:##7|J:##8||{%
\def\SigText{##2}\global\setbox0=\hbox{##2\relax}
\edef\TEST{[\SigText]}
\expandafter\ifx\csname##1\endcsname\TEST\relax\else
\immediate\write16{##1 hat sich geaendert!}\fi
\expandwrite\AUX{\neverexpand\ref{##1}{\TEST}}
\setHI{[\SigText]\ }
\ifnum\HIindent>\maxSpalbr\maxSpalbr\HIindent\fi
\ifnum\HIindent<\altSpalbr\HIindent\altSpalbr\fi
\HI{\marginnote{##1}[\SigText]}
\ifx-##3\relax\else{##3}: \fi
\ifx-##4\relax\else{##4}{\sfcode`.=3000.} \fi
\ifx-##5\relax\else{(##5)} \fi
\ifx-##7\relax\else{##7:} \fi
\ifx-##6\relax\else{##6}\fi
\ifx-##8\relax\else{ ##8}\fi\Par}
\def\Pr|Abk:##1|Sig:##2|Au:##3|Artikel:##4|Titel:##5|Hgr:##6|Reihe:{%
\def\SigText{##2}\global\setbox0=\hbox{##2\relax}
\edef\TEST{[\SigText]}
\expandafter\ifx\csname##1\endcsname\TEST\relax\else
\immediate\write16{##1 hat sich geaendert!}\fi
\expandwrite\AUX{\neverexpand\ref{##1}{\TEST}}
\setHI{[\SigText]\ }
\ifnum\HIindent>\maxSpalbr\maxSpalbr\HIindent\fi
\ifnum\HIindent<\altSpalbr\HIindent\altSpalbr\fi
\HI{\marginnote{##1}[\SigText]}
\ifx-##3\relax\else{##3}: \fi
\ifx-##4\relax\else{##4}{\sfcode`.=3000.} \fi
\ifx-##5\relax\else{In: \it ##5}. \fi
\ifx-##6\relax\else{(##6)} \fi\PrII}
\def\PrII##1|Bd:##2|Verlag:##3|Ort:##4|S:##5|J:##6||{%
\ifx-##1\relax\else{##1} \fi
\ifx-##2\relax\else{\bf ##2}, \fi
\ifx-##4\relax\else{##4:} \fi
\ifx-##3\relax\else{##3} \fi
\ifx-##6\relax\else{##6}\fi
\ifx-##5\relax\else{, ##5}\fi\Par}
\bgroup
\baselineskip12pt
\parskip2.5pt plus 1pt
\hyphenation{Hei-del-berg}
\sfcode`.=1000
\beginsection References. References

}
\def\endrefs{%
\expandwrite\AUX{\neverexpand\ref{Spaltenbreite}{\the\maxSpalbr}}
\ifnum\maxSpalbr=\altSpalbr\relax\else
\immediate\write16{Spaltenbreite hat sich geaendert!}\fi
\egroup}


\marginalnotestrue
\marginalnotesfalse

\def\rhoq{{\overline\rho}}
\def\da{\downarrow}
\def\1{{\textstyle{s\over2}}}
\def\cBq{{\overline\cB}}
\def\cPq{{\overline\cP}}
\def\9{\mathpalette\getstyle}
\def\getstyle#1#2{\setbox0\hbox{$#1#2$}{\vrule width0pt height\ht0}^t\!#2}
\def\tD{{\9D}}

\title{Combinatorics and invariant differential operators on
multiplicity free spaces}

\beginsection Introduction. Introduction

Let $G$ be a connected reductive group acting on a finite dimensional
vector space $U$ (everything defined over $\CC$). We assume that $U$ is a
multiplicity free space, i.e., every simple $G$\_module appears in
$\cP(U)$, the algebra of polynomial functions on $U$, at most
once. Thus, as a $G$\_module,
$$\eqno{}
\cP(U)\cong\bigoplus_{\lambda\in\Lambda_+}M_\lambda
$$
where $\Lambda_+$ is a set of dominant weights and $M_\lambda$ is a
simple $G$\_module of {\it lowest\/} weight $-\lambda$. All elements
of $M_\lambda$ are homogeneous of the same degree, denoted
$\ell(\lambda)$.

Now consider an invariant differential operator $D$ on $U$. It will
act on each irreducible constituent $M_\lambda$ as a scalar, denoted
by $c_D(\lambda)$. It can be shown that $c_D$ extends to a polynomial
function to $V$, the $\CC$\_span of $\Lambda_+$. Thus, $D\mapsto c_D$
is a homomorphism from $\cP\cD(U)^G$, the algebra of invariant
differential operators, into $\cP(V)$. It is possible to determine the
image of this map. One can show, \cite{Montreal}, that there is a
``shift vector'' $\rho\in V$ and a finite reflection group $W\subseteq
GL(V)$ such that the following is an isomorphism:
$$\eqno{E50}
\cP\cD(U)^G\pf\sim\cP(V)^W:D\mapsto p_D(z):=c_D(z-\rho).
$$
Thus, the eigenvalues of $D$ in $\cP(U)$ are the values
$p_D(\rho+\lambda)$, $\lambda\in\Lambda_+$.

The identification \cite{E50} works actually in the much wider context
of $G$\_varieties (see \cite{Annals}) but only multiplicity free
spaces have the following important feature: $\cP\cD(U)^G$ has a
distinguished basis $D_\lambda$, $\lambda\in\Lambda_+$. The
construction of the $D_\lambda$ goes back to Capelli. Via the
identification \cite{E50}, we get also a distinguished basis
$p_\lambda:=p_{D_\lambda}$ of $\cP(V)^W$.

It is possible to characterize the elements of this basis purely in
terms of $V$ without any reference to $U$. Namely,
$p_\lambda$ is the unique $W$\_invariant polynomial on $V$ of degree
$\ell(\lambda)$ which has the following interpolation property:
$$\eqno{}
p_\lambda(\rho+\mu)=\delta_{\lambda\mu}\quad\hbox{for all
}\mu\in\Lambda_+\hbox{ with }\ell(\mu)\le\ell(\lambda).
$$

Note that this is a purely combinatorial description of $p_\lambda$:
all we need to know are $V$, $W$, $\Lambda_+$, $\ell$, and $\rho$.
The first four of these data are rather rigid but there is some
flexibility for $\rho$. In fact, there are many, quite different,
examples of multiplicity free spaces for which $V$, $W$, $\Lambda_+$,
and $\ell$ are the same but $\rho$ is different. This is a
motivation for using the characterization above to {\it define} a
family of polynomials $p_\lambda(z;\rho)$ for an (almost) {\it
arbitrary} $\rho\in V$ (a suggestion of Sahi, see \cite{Sahi}).

In general, not much can be said about $p_\lambda(z;\rho)$ but we
showed in \cite{Existence}\footnote{In fact, the present paper is as a
continuation of \cite{Existence}. For the convenience of the reader we
recalled all relevant results in section~\cite{Capelli}.} that for
$\rho$ in a certain {\it non\_zero\/} subspace $V_0$ of $V$ these
polynomials have remarkable properties. The most important one is the
existence of difference operators $D_h$, $h\in\cP(V)^W$ for which all
polynomials $p_\lambda=p_\lambda(z;\rho)$ are eigenfunctions. More
precisely,
$$\eqno{}
D_h(p_\lambda)=h(\rho+\lambda)p_\lambda.
$$
Thus we can think of the polynomials $p_\lambda(z;\rho)$, $\rho\in
V_0$, as a good deformation of the spectral polynomials
$p_\lambda(z)$.

The central result of the present paper is the Transposition Formula
for $p_\lambda(z;\rho)$. Again, it originates from differential
operators. There, ``transposition'' is the unique {\it
anti\/}automorphism $D\mapsto\tD$ of $\cP\cD(U)$ with
$$\eqno{}
\9x_i=x_i\quad\hbox{and}\quad
\9{\!\!\left({\partial\over\partial x_i}\right)}=-{\partial\over\partial x_i}.
$$
Transposition commutes with the $G$\_action and induces an
automorphism of $\cP\cD(U)^G$. It is a natural problem to calculate
its effect on $\cP(V)^W$ under the identification
\cite{E50}. This is done in section~\cite{Fourier} and the answer is
simply the map $h\mapsto h^-$ where $h^-(z):=h(-z)$
(\cite{Fourier-transform}).

{}From now on, we denote $\cP(V)$ simply by $\cP$. In
section~\cite{Duality}, we compute the action of $h\mapsto h^-$ on
$\cP^W$ with respect to the $p_\lambda$\_basis. The result is the {\it
transposition formula\/} (\cite{Dual2}):
$$\eqno{E61}
q_\lambda(-z)=\sum_\mu (-1)^{\ell(\mu)}p_\mu(\rho+\lambda)q_\mu(z)
\quad\hbox{for all}\quad\lambda\in\Lambda_+.
$$
Here, we used the renormalized polynomials
$$\eqno{}
q_\lambda(z;\rho):={1\over p_\lambda(-\rho;\rho)}\ p_\lambda(z;\rho).
$$
Its proof uses the difference operators $D_h$, an idea which goes
back to Okounkov, \cite{Ok}\footnote{There, the ``transposition
formula'' is called ``binomial formula''.}, who proved it for
shifted Jack polynomials.

A first consequence of the transposition formula is the {\it
evaluation formula\/} (\cite{EvalFormula})
$$\eqno{}
p_\lambda(-\rho;\rho)=(-1)^{\ell(\lambda)}d_\lambda
$$
where
$$\eqno{}
d_\lambda=\prod_{\alpha\in\Delta^+}{\alpha(\rho+\lambda)\over\alpha(\rho)}\ 
\prod_{\omega\in\Phi^+}{(\omega(\rho)+k_\omega)_{\omega(\lambda)}
\over(\omega(\rho)-k_\omega+1)_{\omega(\lambda)}}.
$$
Here $(a)_n=a(a+1)\ldots(a+n-1)$ is the Pochhammer symbol, $\Delta^+$
and $\Phi^+$ are certain finite sets of linear functions on $V$
(positive roots and pseudoroots, respectively), and $k_\omega$ is a
multiplicity function determined by $\rho$.  The number $d_\lambda$ is
called the {\it virtual dimension\/} since, in the case when $\rho$
comes from a multiplicity free space $U$, it computes the dimension of
the irreducible $G$\_modules occurring in $\cP(U)$
(\cite{virdim}). This result generalizes a formula of Upmeier,
\cite{Up}, who considered multiplicity free spaces coming from
Hermitian symmetric spaces (see below).

As already observed in \cite{Ok}, another consequence of the
transposition formula is the {\it interpolation formula}. It gives the
expansion of an arbitrary polynomial $h\in\cP^W$ in terms of the
$p_\lambda$'s. More precisely, we show in section~\cite{Interpolation}
(\cite{Interpol}):
$$\eqno{}
h(z)=\sum_{\lambda\in\Lambda_+}(-1)^{\ell(\lambda)}\widehat h(\rho+\lambda)
p_\lambda(z)
$$
where
$$\eqno{E60}
\widehat h(\rho+\lambda):=\sum_{\mu\in\Lambda_+}(-1)^{\ell(\mu)}
p_\mu(\rho+\lambda)h(\rho+\mu).
$$

Another consequence (also noticed in \cite{Ok}) of the transposition
formula is the symmetry
$$\eqno{E62}
q_\lambda(-\rho-\nu)=q_\nu(-\rho-\lambda),\qquad \lambda,\nu\in\Lambda_+
$$
(just substitute $z=\rho+\nu$ in \cite{E61}). In
section~\cite{Product}, we define a scalar product
on $\cP^W$ by
$$\eqno{}
\<p_\lambda,p_\mu\>:=d_\lambda\delta_{\lambda\mu}.
$$
Then \cite{E62} is explained by the fact
$\<q_\lambda^-,q_\nu^-\>=q_\lambda(-\rho-\nu)$.

Let $\cA\subseteq\|End|_\CC\cP^W$ be the algebra generated by all
multiplication operators $h\in\cP^W$ and all difference operators
$D_h$, $h\in\cP^W$. Then the transformation \cite{E60} can be used to
define an involutory automorphism $X\pfeil\widehat X$ of $\cA$ which
interchanges $h$ and $D_h$ (\cite{InvAut}).  Moreover, we show that
$\cA$ is stable under taking adjoints for the auxiliary
scalar product
$$
\<f,g\>^-:=\<f^-,g^-\>.
$$
More precisely, $(h, D_{h^-})$ is an adjoint pair (\cite{Adjoint1}).

These results are extended in section~\cite{sl-two}. For every
operator $X$ define $X^-$ by $X^-(h):=X(h^-)^-$. Let $\cB$ be the
algebra generated by all $h$, $D_h$, and $D_h^-$ with $h\in\cP^W$. In
other words, $\cB$ is generated by $\cA$ and $\cA^-$. First we observe
that $\cB$ is stable under taking adjoints $X^*$ with respect to the
original scalar product (\cite{adjoints}). The main result of
section~\cite{sl-two} is the construction of a $PGL_2(\CC)$\_action
on $\cB$ which incorporates the two automorphisms $X\mapsto\widehat X$
and $X\mapsto X^-$. For this, let
$$\eqno{}
L:=\ell-D_\ell.
$$
Then we show that $(L,2\ell,L^-)$ forms an $sl_2$\_triple
(\cite{triple}). The $PGL_2(\CC)$\_action is then obtained by integrating the
adjoint action of this triple (\cite{integration}).

In section~\cite{differential}, we study the effect of operators in
$\cB$ on the top homogeneous components of polynomials. More
precisely, both $\cP^W$ and $\cB$ are filtered by degree. Denote their
associated graded algebras by $\cPq^W$ and $\cBq$, respectively. Then
the $\cBq$\_module $\cPq^W$ is called the {\it differential limit\/}
of the $\cB$\_module $\cP^W$.  While it is clear that
$\cPq^W\cong\cP^W$ (since $\cP$ is graded to begin with) we show that
also $\cBq\cong\cB$ (\cite{difflim}). Therefore, the algebra $\cB$ of
{\it difference\/} operators can be replaced by $\cBq$, an algebra of
{\it differential\/} operators (\cite{diffop}). Unfortunately, so far
it seems to be very hard to construct $\cBq$ directly.

In section~\cite{Binomial}, we study another limit, namely the
infinitesimal neighborhood of a particular $W$\_fixed point $\delta$ in
$V$. The transposition formula~\cite{E61} then becomes, in the limit,
a {\it binomial formula\/} (\cite{BinomialTheorem}):
$$\eqno{}
\qq^{(\delta)}_\lambda(z+\delta)=\sum
\limits_{\mu\in\Lambda_+\atop\ell^\delta(\mu)=\ell^\delta(\lambda)}
p_\mu(\rho+\lambda)\ \qq^{(\delta)}_\mu(z)
$$
where $\qq^{(\delta)}_\lambda(z)$ is a certain renormalization of the
top homogeneous component of $p_\lambda(z)$.

Multiplicity free spaces have been classified by Kac \cite{Kac},
Benson\_Ratcliff \cite{BenRat1}, and Leahy \cite{Leahy}. So far,
basically only two classes have been studied in more detail. The case
which got by far the most attention is the so\_called classical
case. It includes the spaces when $G$ is the complexification of the
isotropy group of a Hermitian symmetric space and $U$ is the
complexification of ``$\fp^+$''. Here, the polynomials
$p_\lambda(z;\rho)$ are called {\it shifted Jack polynomials\/} since
their top homogeneous components are the Jack polynomials. By now
there is a rich literature on these polynomials, and most
results of this paper have been previously obtained in that case,
(see, e.g., \cite{Las}, \cite{La2}, \cite{SymCap}, \cite{Ok},
\cite{OO1}, \cite{OlOk}, \cite{Sahi}, \cite{Sek}) even though the
results of section~\cite{sl-two} on the $PGL_2(\CC)$\_action seem to
be new even in the classical case.

The other case, in which the present theory is (mostly) worked out is
the semiclassical case, \cite{Semi}. This includes, e.g., the action
of $GL_n(\CC)$ on $\wedge^2\CC^{n+1}$. Among the few papers which deal
with general multiplicity free spaces are most notably \cite{HoUm},
\cite{Yan}, and \cite{BenRat}.

We found it useful to illustrate most of our results with the case
$\|dim|V=1$. This case is pretty elementary but still quite
interesting. We could have sprinkled specializations to this case all
over the paper but found it more useful to gather everything in a
separate section at the end of the paper. It is recommended to consult
this section frequently in the course of reading this paper or to even
start with it.

\beginsection Fourier. Transposition of differential operators
on multiplicity free spaces

Let $U$ be a finite dimensional complex vector space. Then the algebra
$\cP\cD(U)$ of linear differential operators with polynomial
coefficients has the following presentation: it is generated by $U$
(the directional derivatives) and $U^*$ (the linear functions) which
satisfy the following relations:
$$\eqno{}
[\partial_1,\partial_2]=0, [x_1,x_2]=0, [\partial,x]=\partial(x)\quad
\hbox{for all }\partial_1,\partial_2,\partial\in U,\quad x_1,x_2,x\in U^*.
$$
This implies that there is a unique {\it anti\/}automorphism
$D\mapsto\tD$ of $\cP\cD(U)$ with
$$\eqno{}
\9\partial=-\partial,\quad
\9x=x\quad\hbox{for all $u\in U$ and $x\in U^*$}.
$$
The operator $\tD$ is called the {\it transpose\/} of $D$.

Let $G$ be an algebraic group $G$ acting linearly on $U$. Then transposition
is $G$\_equivariant. It follows, that it induces an
antiautomorphism of the algebra $\cP\cD(U)^G$ of $G$\_invariant
differential operators.

Now assume that $G$ is connected, reductive and $U$ a multiplicity
free space. This means that the algebra $\cP(U)$ of regular functions
is multiplicity free as a $G$\_module. Then it is easy to show that
$\cP\cD(U)^G$ is commutative (in fact, this will be shown
below). Thus, transposition is an automorphism of $\cP\cD(U)^G$. The
purpose of this section is to calculate this automorphism explicitly.

To do this, we need first an explicit description of the algebra
$\cP\cD(U)^G$ itself. Fix a Borel subgroup $B$ of $G$ and a maximal torus
$T$ of $B$. By assumption, the algebra $\cP(U)$ decomposes as a
$G$\_module as $\cP(U)=\oplus_{\lambda\in\Lambda_+}\cP^\lambda$ where
$\Lambda_+\subseteq(\|Lie|T)^*$ is a certain set of integral dominant
weights and $\cP^\lambda$ is a simple $G$\_module with {\it lowest\/}
weight $-\lambda$.

Every $D\in\cP\cD(U)^G$ acts on $\cP^\lambda$ as multiplication by a
scalar, which is denoted $c_D(\lambda)$. Let $V\subseteq(\|Lie|T)^*$
be the $\CC$\_span of $\Lambda_+$. Then $c_D$ is the restriction of a
unique polynomial function on $V$ (also denoted by $c_D$) to
$\Lambda_+$ (see \cite{Montreal}~Cor. 4.4). Thus, we obtain an embedding
$$\eqno{}
\cP\cD(U)^G\into\cP(V):D\mapsto c_D
$$
which shows, in particular, that $\cP\cD(U)^G$ is commutative.

To describe the image of this embedding we need some more notation.
Let $P\supseteq B$ be the largest parabolic subgroup such that all
elements of $\Lambda_+$, which are characters of $T$, extend to
characters of $P$. Let $\beta$ be the sum of all roots in the
unipotent radical of $P$ and let $\chi$ be the sum of all weights of
$U$. Then it is shown in \cite{Existence}~\S7 that
$\rho:={1\over2}(\beta+\chi)\in V$. Using this weight, we define a new
embedding
$$\eqno{}
\cP\cD(U)^G\into\cP(V):D\mapsto p_D
$$
where
$$\eqno{E2}
p_D(v):=c_D(v-\rho).
$$
Then we have the following (Harish Chandra) isomorphism:

\Theorem HarishChandra. {\rm(\cite{Montreal} Theorem~4.8)} There is a
unique finite subgroup $W\subseteq GL(V)$ such that $D\mapsto p_D$
establishes an isomorphism between $\cP\cD(U)^G$ and $\cP(V)^W$.

\noindent Now we can make transposition of invariant differential
operators explicit:

\Theorem Fourier-transform. Let $U$ be a multiplicity free space for
$G$. Then
$$\eqno{E1}
p_{\tD}(v)=p_D(-v)\quad\hbox{for every}
\quad D\in\cP\cD(U)^G\ {\rm and}\ v\in V.
$$

\Proof: Using \cite{E2} we have to prove
$$\eqno{E3}
c_\tD(v)=c_D(-v-\chi-\beta).
$$
Let $\fZ(\fg)$ be the center of the universal enveloping algebra
$\fU(\fg)$ of the Lie algebra of $G$. The action of $G$ on $U$ induces
a homomorphism $\Psi:\fU(\fg)\pfeil\cP\cD(U)$ which maps $\fZ(\fg)$ to
$\cP\cD(U)^G$. We are going to verify \cite{E3} first for operators in
the image of $\Psi$.

Let $u_i$ be a basis of $U$ where each $u_i$ is a weight vector with
weight $\chi_i$. Let $x_i\in U^*$ be the dual basis and
$\partial_i:={\partial\over\partial x_i}$. Consider the decomposition
$\fg=\fn\oplus\ft\oplus\fn^-$. For $\eta\in\ft$ we have
$\Psi(\eta)=-\sum_i\chi_i(\eta)x_i\partial_i$. Thus
$$\eqno{}
\9\Psi(\eta)=-\sum_i\chi_i(\eta)(-\partial_i)x_i
=\sum_i\chi_i(\eta)(x_i\partial_i+1)=-\Psi(\eta)+\chi(\eta).
$$
If $\eta\in\fn^\pm$ then $\Psi(\eta)=\sum_{i\ne
j}a_{ij}x_i\partial_j$, hence
$\9\Psi(\eta)=-\Psi(\eta)$. Now observe that $\chi$, being
the sum of all weights of $U$, is actually a character of all of
$\fg$. Thus, we can define an {\it anti\/}automorphism $\tau$ of
$\fU(\fg)$ by $\tau(\eta):=-\eta+\chi(\eta)$ for all $\eta\in\fg$
and the discussion above showed
$$\eqno{}
\9\Psi(\xi)=\Psi(\tau(\xi))\quad\hbox{for all}\quad\xi\in\fU(\fg).
$$

Let $\xi\in\cZ(\fg)$ and $D=\Psi(\xi)$. From the theorem of
Poincar\'e\_Birkhoff\_Witt follows that $\xi$ decomposes uniquely as
$\xi=\xi_0+\xi_1$ with $\xi_0\in\fU(\ft)$ and
$\xi_1\in\fn^-\fU(\fg)\fn$. Since $\fU(\ft)=S(\ft)$, we can regard
$\xi_0$ as a function on $\ft^*$ and write $\xi_0(v)$ for its
value at $v\in\ft^*$. In particular, we have
$$\eqno{}
\tau(\xi_0)(v)=\xi_0(-v+\chi).
$$
On the other hand $\tau(\xi)=\tau(\xi_0)+\tau(\xi_1)$ with
$\tau(\xi_0)\in\fU(\ft)$ and $\tau(\xi_1)\in\fn\fU(\fg)\fn^-$. Let
$f$ be a lowest weight vector of $\cP^\lambda$. By definition, it has
weight $-\lambda$. Thus $\Psi(\tau(\xi_1))f=0$, and we have
$$\eqno{}
\tD f=\Psi(\tau(\xi_0))f=\tau(\xi_0)(-\lambda)f=\xi_0(\lambda+\chi)f.
$$
Thus
$$\eqno{E5}
c_\tD(\lambda)=\xi_0(\lambda+\chi)
\quad\hbox{for all}\quad\lambda\in\Lambda_+.
$$

Let $w_0$ be the longest element of the Weyl group $\Wq$ of $G$. Then
the highest weight of $\cP^\lambda$ is $-w_0\lambda$. Thus
$c_D(\lambda)=\xi_0(-w_0\lambda)$. Let $\rhoq$ be the half\_sum of
positive roots of $G$. By the (original) Harish Chandra isomorphism,
the function $v\mapsto\xi_0(v-\rhoq)$ is $\Wq$\_invariant. From
$w_0\rhoq=-\rhoq$, we get
$$\eqno{E4}
c_D(\lambda)=\xi_0(-w_0\lambda)=
\xi_0(w_0(-\lambda-\rhoq)-\rhoq)=\xi_0(-\lambda-2\rhoq).
$$

Let $L$ be the Levi complement of $P$ and $w_L$ the longest element of
its Weyl group. Then we have the relation
$\beta=\rhoq+w_L\rhoq$. Since $\Lambda_+$ is Zariski dense in $V$,
equation \cite{E4} is valid for all $\lambda\in V$. In particular, we
can replace $\lambda$ by $-\lambda-\chi-\beta$. Thus,
$$\eqno{E6}
c_D(-\lambda-\chi-\beta)=\xi_0(\lambda+\chi+\beta-2\rhoq)=
\xi_0(\lambda+\chi+w_L\rhoq-\rhoq).
$$
Now we use the fact that $\lambda\in V$ and $\chi$ are
$w_L$\_fixed. Hence
$$\eqno{E7}
\xi_0(\lambda+\chi+w_L\rhoq-\rhoq)=\xi_0(w_L(\lambda+\chi+\rhoq)-\rhoq)=
\xi_0(\lambda+\chi).
$$
Equations \cite{E5}, \cite{E6}, and \cite{E7} imply \cite{E3} for
$D=\Psi(\xi)$.

Now we consider the general case. Clearly, there is a unique
automorphism $\sigma$ of $\cP^W$ such that
$$\eqno{}
\sigma(c_D)(v)=c_\tD(-v)
\quad\hbox{for all}\quad D\in\cP\cD(U)^G
$$
and we have to show that $\sigma$ is the identity. By what we proved
above, $\sigma$ fixes the subalgebra $\cP_0:=\{c_D\mid
D\in\Psi\fZ(\fg)\}$ pointwise. Since $\cP(\ft^*)$ is finitely
generated as a $\fZ(\fg)=\cP(V)^\Wq$\_module, also $\cP(V)^W$ is a
finitely generated $\cP_0$\_module. Let $\cK$ be the quotient field of
$\cP(V)^W$. Then we see that $[\cK:\cK^{\<\sigma\>}]$ is finite which
implies that $\sigma$ has finite order.

In the last step, we use that transposition is filtration
preserving. More precisely, $\cP\cD(U)^G$ is filtered by the order of
a differential operator and $\cP(V)^W$ is filtered by degree. The
associated graded ring of $\cP\cD(U)^G$ is $\cP(U\oplus U^*)^G$ and
transposition induces on the latter the action $(u,\alpha)\mapsto
(u,-\alpha)$. The map $D\mapsto c_D$ is degree preserving. Thus
transposition acts on $\|gr|\cP(V)^W$ by $v\mapsto-v$. This shows that
$\sigma$ acts on $\|gr|\cP(V)^W$ as identity. But $\sigma$ has finite
order, hence is linearly reductive. This implies that $\sigma$ is the
identity on $\cP(V)^W$.\qed

\beginsection Capelli. Capelli polynomials

This section is a synopsis of the essential parts of
\cite{Existence}. We have seen that to every multiplicity free space
there is attached a finite dimensional vector space $V$, a finite
reflection group $W$ acting on it and a finitely generated monoid
$\Lambda_+$ of dominant weights. Additionally, we have a linear
function $\ell:V\pfeil\CC$ such that $\ell(\lambda)=\|deg|f$ for any
non\_zero $f\in\cP^\lambda$. These data are by no means unrelated and
in \cite{Existence} we proposed a set of axioms which we are not going
to repeat since we rarely need them directly. From now on we forget
about multiplicity free spaces and consider just
structures\footnote{Actually, in \cite{Existence} we found it more
convenient to state the axioms in terms of the equivalent data
$(\Gamma,\Sigma,W,\ell)$.} $(V,W,\Lambda_+,\ell)$ satisfying these
axioms. Note that all multiplicity free actions are classified,
\cite{Kac}, \cite{BenRat1}, \cite{Leahy}. The ensuing combinatorial
structures are described in \cite{Existence}~\S8.

\def\dsseq{\hbox{$\cap$\vrule}}
Inside $V$ we are going to consider the following objects:
$$\eqno{}
\matrix{\Sigma^\vee&&\Lambda_1\cr
\dsseq&&\dsseq\cr
\Lambda_+&\subseteq&\Lambda&\subseteq&\Gamma^\vee&\subseteq&V\cr}
$$
Here $\Gamma^\vee$ is the lattice generated by $\Lambda_+$ and
$\Lambda$ is the submonoid generated by all $w\eta$ with $w\in W$ and
$\eta\in\Lambda_+$. The minimal set of generators of $\Lambda_+$ is
denoted by $\Sigma^\vee$. It forms a basis of $\Gamma^\vee$ and
$V$. Also $\Lambda$ has a minimal set of generators which is denoted
by $\Lambda_1$. It coincides with the set of all $w\eta$ with $w\in
W$, $\eta\in\Sigma^\vee$ and $\ell(\eta)=1$.

Inside the dual space $V^\vee$ we need the following objects:
$$\eqno{}
\matrix{\Sigma&&\Delta\cr
\dsseq&&\dsseq\cr
\Phi&\subseteq&\Gamma&\subseteq& V^\vee\cr}
$$

Here $\Gamma$ is the lattice dual to $\Gamma^\vee$ and $\Sigma$ is the
dual basis of $\Sigma^\vee$. The elements of $\Phi:=\cup_{w\in
W}w\Sigma$ are called {\it pseudoroots.} Attached to the reflection
group $W$ there is a unique root system $\Delta$ such that all roots
are primitive vectors.

Let $\pm W$ be the group generated by $W$ and $-1$. Then we define
$$\eqno{}
V_0:=\{\rho\in V\mid\hbox{For all $\omega_1,\omega_2\in\Sigma$
with $\omega_1\in \pm W\omega_2$ holds $\omega_1(\rho)=\omega_2(\rho)$}\}
$$
Thus, for $\rho\in V_0$ and for every $\omega\in\Phi\cup(-\Phi)$ we can define
$k_\omega:=\omega_1(\rho)$ where $\omega_1\in
\pm W\omega\cap\Sigma$. In particular we have $k_\omega=k_{-\omega}$
for all $\omega\in\Phi$.

\Examples: For the rank one case see section \cite{Beispiel}. Here we
illustrate the notation above in two examples. In
the {\it classical case\/} we have:

\noindent $V:=\CC^n$, $W:=S_n$ (symmetric
group),$\Lambda_+:=\{\lambda\in\ZZ^n\mid\lambda_1\ge\ldots\ge\lambda_n\ge0\}$
(partitions), and $\ell(\lambda):=\sum_i\lambda_i$.
Let $e_1,\ldots,e_n$ be the canonical basis of $\CC^n$ and
$z_1,\ldots,z_n\in(\CC^n)^\vee$ its dual basis. Then
we get the following derived data:
\medskip
\halign{$#=$\hfill&$#$\hfill\qquad&$#=$\hfill&$#$\hfill\cr
\multispan2Subsets of $V$\hfill&\multispan2Subsets of $V^\vee$\hfill\cr
\noalign{\vskip-7pt}
\multispan2\hrulefill\qquad&\multispan2\hrulefill\cr
\Gamma^\vee&\ZZ^n&\Gamma&\ZZ^n\cr
\Sigma^\vee&\{e_1+\ldots+e_i\mid 1\le i\le n\}&
\Sigma&\{z_i-z_{i+1}\mid 1\le i<n\}\cup\{z_n\}\cr
\Lambda_1&\{e_i\mid 1\le i\le n\}&\Phi&
\{z_i-z_j\mid 1\le i\ne j\le n\}\cup\{z_i\mid 1\le i\le n\}\cr
\Lambda&\NN^n&\Delta&
\{z_i-z_j\mid 1\le i\ne j\le n\}\cr
V_0&\{\sum_i[(n-i)r+s]e_i\mid r,s\in\CC\}&\multispan2\hfill\cr}
\medskip\noindent
Observe that $\Delta$ is a subset of $\Phi$. This makes the classical
case rather exceptional. It has been the topic of the
papers \cite{SymCap} and \cite{OlOk} among others.

\medskip\goodbreak
The second example is the {\it semiclassical case}. Then:

\noindent
$V:=\CC^n$, $W:=\{w\in S_n\mid w(i)-i\hbox{ even for all }i\}$
(semisymmetric group),\hfill\break
$\Lambda_+:=\{\lambda\in\ZZ^n\mid\lambda_1\ge\ldots\ge\lambda_n\ge0\}$
(partitions), and $\ell(\lambda):=\sum_{i\ \rm odd}\lambda_i$.

\noindent We get the following derived data:
\medskip

\halign{$#=$\hfill&$#$\hfill\qquad&$#=$\hfill&$#$\hfill\cr
\multispan2Subsets of $V$\hfill&\multispan2Subsets of $V^\vee$\hfill\cr
\noalign{\vskip-7pt}
\multispan2\hrulefill\qquad&\multispan2\hrulefill\cr
\Gamma^\vee&\ZZ^n&\Gamma&\ZZ^n\cr
\Sigma^\vee&\{e_1+\ldots+e_i\mid 1\le i\le n\}&
\Sigma&\{z_i-z_{i+1}\mid 1\le i<n\}\cup\{z_n\}\cr
\Lambda_1&\{e_i\mid i\hbox{ odd}\}\cup\{e_i{+}e_j\mid i\hbox{ odd},
j\hbox{ even}\}&\Phi&\{z_i-z_j\mid i-j\ {\rm odd}\}\cup
\{z_i\mid n-i\ {\rm even}\}\cr
\Lambda&\{\lambda\in\NN^n\mid
\sum_{i\ \rm odd}\lambda_i\ge\sum_{i\ \rm even}\lambda_i\}&\Delta&\{z_i-z_j\mid i\ne j, i-j\ {\rm even}\}\cr
V_0&\{\sum_i[(n-i)r+s]e_i\mid r,s\in\CC\}&\multispan2\hfill\cr}

\medskip\noindent The semiclassical case has been investigated in \cite{Semi}.

\medskip
We are going to need the following non\_degeneracy conditions for
$\rho\in V_0$. Let
$\Delta^+:=\{\alpha\in\Delta\mid\alpha(\Sigma^\vee)\ge0\}$. Then
$$\eqno{}
\hbox{$\rho\in V_0$ is
$\left\{\vcenter{\baselineskip12pt\hbox{\it dominant}
                 \hbox{\it non-integral}}
\right\}$ if  $\alpha(\rho)\not\in
\left\{\vcenter{\baselineskip12pt\hbox{$\ZZ_{<0}$}
                \hbox{$\ZZ$}}
\right\}$ for all $\alpha\in\Delta^+$.}
$$

Let $\cP$ denote the algebra of polynomial functions on $V$.  The next
theorem introduces one of the main objects of the theory: a
distinguished basis of $\cP^W$ whose elements are sometimes called
Capelli polynomials since they are related to the Capelli identities.

\Theorem Main-Prop-E. {\rm (\cite{Existence}~Theorem 3.6)} Let $\rho\in
V_0$ be dominant.
\item{a)} For every $\lambda\in\Lambda_+$ there is a unique polynomial
$p_\lambda\in\cP^W$ with $\|deg|p_\lambda\le\ell(\lambda)$ and
$p_\lambda(\rho+\mu)=\delta_{\lambda\mu}$ (Kronecker delta) for all
$\mu\in\Lambda_+$ with $\ell(\mu)\le\ell(\lambda)$.
\item{b)}For every $d\in\NN$, the set of $p_\lambda$ with
$\ell(\lambda)\le d$ forms a basis of the space of $p\in\cP^W$ with
$\|deg|p\le d$.\Par

\noindent The polynomials vanish, in fact, in many more points than
they are supposed to. This is the content of the Extra Vanishing
Theorem:

\Theorem ExtraVanishing. {\rm(\cite{Existence} Corollary~3.9)} Let
$\rho\in V_0$ be dominant. Then for any $\lambda,\mu\in\Lambda_+$
holds $p_\lambda(\rho+\mu)=0$ unless $\mu\in\lambda+\Lambda$.

For $d\in\ZZ$ define the following variant of the falling factorial
polynomial:
$$\eqno{}
[z\da d]:=\cases{z(z-1)\ldots(z-d+1)&if $d>0$,\cr
1&otherwise.\cr}
$$
Then, for every $\tau\in\Gamma$ we define the
rational function
$$\eqno{Def-f-l}
f_\tau(z):={
\prod\limits_{\omega\in\Phi}[\omega(z)-k_\omega\da\omega(\tau)]
\over
\prod\limits_{\alpha\in\Delta}[\alpha(z)\da\alpha(\tau)]}.
$$
One of its main features are the following cut\_off properties:

\Lemma Cut-Off. {\rm(\cite{Existence} Lemmas 3.2 and 5.2)} Let $\rho\in
V_0$ be non\_integral and $\tau\in\Lambda$.
\item{a)}Assume $\lambda\in\Lambda_+$ but
$\mu:=\lambda-\tau\not\in\Lambda_+$. Then $f_\tau(\rho+\lambda)=0$.
\item{b)}Assume $\mu\in\Lambda_+$ but
$\lambda:=\mu+\tau\not\in\Lambda_+$. Then $f_\tau(-\rho-\mu)=0$.\Par

For any $\eta\in V$ we define the shift operator $T_\eta$ on $\cP$
by $(T_\eta f)(z)=f(z-\eta)$. Then the difference operator
$$
L:=\sum_{\eta\in\Lambda_1}f_\eta(z)T_\eta
$$
has very remarkable properties. Since its coefficients are rational
functions, it doesn't act on $\cP$ but is does on $\cP^W$.

\Examples: 1. Classical case:
$$\eqno{}
L=\sum_{i=1}^n\left[\prod_{j\ne i}{z_i-z_j-r\over z_i-z_j}\right]
(z_i-s)\ T_{e_i}.
$$

\noindent 2. Semiclassical case:
$$\eqno{}
\eqalign{L=&\sum_{i\ \rm odd}\left[{\prod_{j\ \rm even}(z_i-z_j-r)\over
\prod_{j\ne i\ \rm odd}(z_i-z_j)}\right](z_i-s)\ T_{e_i}+\cr
&+\sum_{i\ {\rm odd}\atop j\ {\rm even}}\left[
{\prod_{k\ne j\ \rm even}(z_i-z_k-r)\prod_{k\ne i\ \rm odd}(z_j-z_k-r)
\over
\prod_{k\ne i\ \rm odd}(z_i-z_k)\prod_{k\ne j\ \rm even}(z_j-z_k)}\right]
(z_i-s)\ T_{e_i+e_j}.\cr}
$$
\medskip\noindent
One of the main properties of $L$ is:

\Theorem nilp. {\rm(\cite{Existence}~Corollary~5.7)} Consider
$h\in\cP^W$ as multiplication operator on $\cP^W$. Then
$(\|ad|L)^n(h)=0$ for $n>\|deg|h$.

\noindent Thus, for every $h\in\cP^W$ we can define the difference
operator
$$\eqno{}
D_h:=\|exp|(\|ad|L)(h).
$$
The most important special case is the {\it difference Euler
operator\/} $E:=D_\ell=\ell-L$. All these operators are diagonalized
by the $p_\lambda$. More precisely:

\Theorem. {\rm(\cite{Existence} Theorem~5.8)} Let $h\in\cP^W$. Then
$$\eqno{}
D_h(p_\lambda)=h(\rho+\lambda)p_\lambda
\quad\hbox{for all}\quad\lambda\in\Lambda_+.
$$

\noindent
In the classical and semiclassical case, these difference operators
have been first constructed explicitly in \cite{SymCap} and
\cite{Semi}, respectively. In general, much less is known. The rough
structure of $D_h$ is explained by the following Lemma.

\Lemma roughstructure. There is an expansion
$$\eqno{}
D_h=\sum_\eta b_\eta^h(z)T_\eta
$$
where $b_\eta^h(z)$ is a rational function and $\eta\in\Lambda$ with
$\ell(\eta)\le\|deg|h$.

\Proof: That $b_\eta^h(z)$ is obvious from the definition. Let
$d=\|deg|h$. Then
$$\eqno{}
D_h=\sum_\eta b_\eta^h(z)T_\eta=\sum\limits_{n=0}^d{1\over n!}(\|ad|L)^n(h)
$$
by \cite{nilp}. Thus $b_\eta^h=0$ unless $\eta$ is the sum of at most $d$
elements of $\Lambda_1$. But this implies $\eta\in\lambda$ with
$\ell(\eta)\le d$.\qed

\noindent There is a strong connection between the difference
operators $D_h$ and Pieri\_type formulas. For this we define for every
$\lambda\in\Lambda_+$ the {\it virtual dimension}\footnote{See
\cite{virdim} for an explanation of this term.} as
$$\eqno{E30}
d_\lambda:=(-1)^{\ell(\lambda)}
{f_\lambda(-\rho)\over f_\lambda(\rho+\lambda)}.
$$
It can be rewritten as
$$\eqno{E51}
d_\lambda=\prod_{\alpha\in\Delta^+}{\alpha(\rho+\lambda)\over\alpha(\rho)}\ 
\prod_{\omega\in\Phi^+}{(\omega(\rho)+k_\omega)_{\omega(\lambda)}
\over(\omega(\rho)-k_\omega+1)_{\omega(\lambda)}}.
$$
where $\Phi^+:=\{\omega\in\Phi\mid\omega(\Sigma^\vee)\ge0\}$. Thus,
the following condition on $\rho$ is designed to make sure that
$d_\lambda$ is defined and non\_zero: we call $\rho$ {\it strongly
dominant\/} if for all $\alpha\in\Delta^+$ and $\omega\in\Phi^+$:
$$\eqno{}
\alpha(\rho)\not\in\ZZ_{\le0},\quad
\omega(\rho)-k_\omega\not\in\ZZ_{<0},\quad
\omega(\rho)+k_\omega\not\in\ZZ_{\le0}.
$$

\Remark: All $\rho$'s coming from multiplicity free actions are
strongly dominant.

\Theorem Pieri1. Let $\rho$ be strongly dominant and
non\_integral. Let $h\in\cP^W$. Then
$$\eqno{P-Formula-1}
h(-z)p_\mu(z)=\sum_\tau (-1)^{\ell(\tau)}{d_\mu\over d_{\mu+\tau}}
b_\tau^h(-\rho-\mu)p_{\mu+\tau}(z)\quad\hbox{for every}\quad\mu\in\Lambda_+.
$$
Here, the sum runs over those $\tau\in\Lambda$ with $\mu+\tau\in\Lambda_+$.

\Proof: This is the combination of formulas (7), (8), and (13) of
\cite{Existence}.\qed

\noindent Later, we are also going to need the following more explicit
Pieri type formula..

\Theorem. {\rm(\cite{Existence}~Corollary~3.11)}
Let $\lambda\in\Lambda_+$ and $k\in\NN$. Then
$$\eqno{Pieriell}
{\ell(z)-\ell(\rho+\lambda)\choose k}p_\lambda(z)=
\sum_{\mu\in\Lambda_+\atop\ell(\mu-\lambda)=k}p_\lambda(\rho+\mu)p_\mu(z).
$$

\beginsection Duality. The transposition formula

In section~\cite{Fourier}, we showed the representation theoretic meaning of
the transformation $h(z)\mapsto h(-z)$ on $\cP^W$. Now, we would like
to express it in terms of the basis $p_\lambda$.

Difference operators act naturally on function from the {\it
left}. Now, we consider also their action on (finite) measures on the
{\it right}. More precisely, for any $v\in V$ let
$\delta_v:\cP\pfeil\CC:f\mapsto f(v)$ be the evaluation map
(a.k.a. Dirac measure). Then the difference operator $D=\sum_\eta
a_\eta(z) T_\eta$ acts on $\delta_v$ by
$$\eqno{}
\delta_vD:=\sum_\eta a_\eta(v)\delta_{v-\eta},
$$
provided the coefficient functions $a_\eta$ are defined in
$z=v$. In that case, we have $(\delta_vD)(f)=\delta_v(D(f))$. We are
interested in measures supported in points of the form $-\rho-\mu$,
$\mu\in\Lambda_+$. Therefore, we define for every $d\in\NN$ the space
$$\eqno{}
\Mq_d:=\bigoplus\limits_{\mu\in\Lambda_+\atop\ell(\mu)\le
d}\CC\,\delta_{-\rho-\mu}\quad\hbox{and}\quad\Mq:=\cup_d\Mq_d
$$

\Proposition FiltIso. Let $\rho$ be strongly dominant and
non\_integral. Then $\Mq$ is $D_h$\_stable for all
$h\in\cP^W$. Moreover, the map
$$\eqno{}
\phi:\cP^W\pfeil\Mq:h\mapsto\delta_{-\rho}D_h
$$
is an isomorphism of filtered $\CC$\_vector spaces.

\Proof: The non\_integrality of $\rho$ makes sure that $\delta D_h$ is
defined for every $\delta\in\Mq$. Clearly, the space $\Mq$ is stable
for the multiplication operator $h$. Thus it suffices to show that
$\Mq$ is $L$\_stable. Since $\delta_{-\rho-\mu}L=\sum_\eta
f_\eta(-\rho-\mu)\delta_{-\rho-\mu-\eta}$, we have to show: for every
$\mu\in\Lambda_+$ holds $\mu+\eta\in\Lambda_+$ or
$f_\eta(-\rho-\mu)=0$. But this is a special case of \cite{Cut-Off}{\it b}.

\cite{roughstructure} implies that $\phi$ preserves filtrations. Since
the filtration spaces on both sides are of the same finite dimension
(\cite{Main-Prop-E}{\it b}) it suffices to show that $\phi$ is
injective. If $\phi(h)=0$ then $b_\tau^h(-\rho)=0$ for all
$\tau\in\Lambda_+$. \cite{Pieri1}, applied to $\mu=0$, then implies
$h=0$.\qed

\noindent
The following consequence is needed in the proof of \cite{Dual2}. 
A much stronger result will proved later on (\cite{EvalFormula}).

\Corollary AAA. Let $\rho$ be strongly dominant and
non\_integral\footnote{See \cite{BBB}.}. Then
$p_\lambda(-\rho)\ne0$ for all $\lambda\in\Lambda_+$.

\Proof: Suppose $p_\lambda(-\rho)=0$. Using the bijectivity of
$\phi$ we get for every $\mu\in\Lambda_+$ a function $h\in\cP^W$ with
$\delta_{-\rho}D_h=\delta_{-\rho-\mu}$. Hence
$$\eqno{}
p_\lambda(-\rho-\mu)=\delta_{-\rho-\mu}p_\lambda=
\delta_{-\rho}D_h(p_\lambda)=h(\rho+\lambda)p_\lambda(-\rho)=0.
$$
Since $-\rho-\Lambda_+$ is Zariski dense in $V$ we conclude
$p_\lambda=0$ which is not true.\qed

\noindent
It is convenient to renormalize $p_\lambda$ such that its value at
$-\rho$ becomes $1$. Therefore, put
$$\eqno{}
q_\lambda(z):={p_\lambda(z)\over p_\lambda(-\rho)}.
$$
Then we can formulate the {\it transposition formula}:

\Theorem Dual2. Let $\rho$ be strongly dominant and
non\_integral\footnote{See \cite{BBB}.}. Then
$$\eqno{DualFormula}
q_\lambda(-z)=\sum_\mu (-1)^{\ell(\mu)}p_\mu(\rho+\lambda)q_\mu(z)
\quad\hbox{for all}\quad\lambda\in\Lambda_+.
$$

\Proof: The polynomials $q_\mu(-z)$ form also a basis of $\cP^W$. Thus, every
$f\in\cP^W$ has an expansion
$$\eqno{E9}
f(z)=\sum_{\mu\in\Lambda_+}a_\mu(f)q_\mu(-z)
$$
where $a_\mu$ is a linear function on $\cP^W$. We claim
$a_\mu\in\Mq_{\ell(\mu)}$. To see that we evaluate \cite{E9} in
$z=-\rho-\mu$ and get
$$\eqno{}
\delta_{-\rho-\mu}(f)=\sum_\tau a_\tau(f)q_\tau(\rho+\mu)=
p_\mu(-\rho)^{-1}a_\mu(f)+\sum_{\ell(\tau)<\ell(\mu)}a_\tau(f)q_\tau(\rho+\mu)
$$
The second equation holds by \cite{Main-Prop-E}{\it a)}. Now the claim
follows by induction.

The claim and \cite{FiltIso} imply that for every $\mu\in\Lambda_+$
there is $h_\mu\in\cP^W$ with $\|deg|h_\mu\le\ell(\mu)$ and
$a_\mu(f)=\delta_{\rho}D_{h_\mu}f=(D_{h_\mu}f)(-\rho)$. Applying this
to $f=q_\lambda$ yields
$$\eqno{}
a_\mu(q_\lambda)=(D_\mu q_\lambda)(-\rho)=h_\mu(\rho+\lambda).
$$
On the other hand, $q_\lambda(z)$ and
$(-1)^{\ell(\lambda)}q_\lambda(-z)$ have the same top homogeneous
component. Thus we get directly from \cite{E9} that
$$\eqno{}
a_\mu(q_\lambda)=(-1)^{\ell(\mu)}\delta_{\lambda\mu}\quad\hbox{for
all}\quad\lambda,\mu\in\Lambda_+\ \hbox{with}\ \ell(\lambda)\le\ell(\mu).
$$
Thus $(-1)^{\ell(\mu)}h_\mu$ matches the definition of $p_\mu$ which
implies
$a_\mu(q_\lambda)=(-1)^{\ell(\mu)}p_\mu(\rho+\lambda)$. Inserting this
into \cite{E9} and replacing $z$ by $-z$ gives formula
\cite{DualFormula}.\qed

\Remark: In the classical case, the transposition formula was first
proved by Okounkov in \cite{Ok} and Lassalle \cite{La2} (even in the
Macdonald polynomial setting). There is was called a ``binomial
theorem'' but we prefer to reserve this term to the limiting case
discussed in section~\cite{Binomial}. We followed Okounkov's approach
to the transposition formula with some substantial modifications. In
particular, we don't need to know the difference operators very
explicitly. The semiclassical case was done in \cite{Semi}.

\medskip \noindent A first consequence of the transposition formula is the
following symmetry result:

\Corollary Symm. Let $\rho$ be strongly dominant and
non\_integral\footnote{See \cite{BBB}.}. Then
$$\eqno{Symmetrisch}
q_\lambda(-\rho-\nu)=q_\nu(-\rho-\lambda)\quad\hbox{for
all $\lambda,\nu\in\Lambda_+$.}
$$

\Proof: Evaluate the transposition formula~\cite{DualFormula} in
$z=\rho+\nu$. Then the right-hand side is symmetric in $\lambda$ and
$\nu$.\qed

\noindent From this, we derive a Pieri formula for the $q_\mu$:

\Theorem Pieri2. Let $\rho\in V_0$ be strongly dominant and
non\_integral. Then
$$\eqno{pieriII}
h(-z)q_\mu(z)=\sum_{\tau\in\Lambda}b_\tau^h(-\rho-\mu)q_{\mu+\tau}(z)
\quad\hbox{for every $h\in\cP^W$}.
$$

\Proof: Consider the eigenvalue equation for $D_h$:
$$\eqno{}
\sum_\tau b_\tau^h(z)q_\nu(z-\tau)=D_h(q_\nu)=h(\rho+\nu)q_\nu(z).
$$
Now substitute $z=-\rho-\mu$ and apply
\cite{Symmetrisch} to both sides:
$$\eqno{}
\sum_\tau b_\tau^h(-\rho-\mu)q_{\mu+\tau}(-\rho-\nu)=
h(\rho+\nu)q_\mu(-\rho-\nu)
$$
(if $\mu+\tau\not\in\Lambda_+$ then $b_\tau^h(-\rho-\nu)=0$, see
\cite{Existence} Proposition~6.3).  This implies \cite{pieriII} since
$-\rho-\Lambda_+$ is Zariski dense in $V$.\qed

\noindent By comparing formulas \cite{P-Formula-1} and \cite{pieriII}
we obtain the {\it evaluation formula}:

\Corollary EvalFormula. Let $\rho\in V_0$ be strongly dominant. Then for
all $\mu\in\Lambda_+$ holds
$$\eqno{E8}
p_\mu(-\rho)=(-1)^{\ell(\mu)}d_\mu.
$$

\Proof: Assume first that $\rho$ is non\_integral. We apply the Pieri
formula \cite{P-Formula-1} to $h(z)=p_\lambda(-z)$ and $\mu=0$. Since
$p_0=1$ we get
$$\eqno{E52}
p_\lambda(-z)\cdot1={(-1)^{\ell(\lambda)}\over d_\lambda}
b_\lambda^h(-\rho)p_\lambda(z)+\hbox{lower order terms}.
$$
Doing the same thing with \cite{pieriII} gives
$$\eqno{}
p_\lambda(-z)\cdot1=b_\lambda^h(-\rho)q_\lambda(z)+\hbox{lower order terms}
$$
Comparing these two formulas proves the evaluation formula. It follows
{}from \cite{E51} that both sides of \cite{E8} are defined when $\rho$
is just strongly dominant. Thus, we can drop the non\_integrality
assumption by a continuity argument.\qed

The last argument of the preceding proof gives:

\Corollary BBB. In \cite{AAA}, \cite{Dual2}, and \cite{Symm} it suffices
to assume that $\rho$ is strongly dominant.

\Remark: This refinement is important since $\rho$\_vectors coming
{}from multiplicity free spaces are almost never non\_integral.
\medskip
\noindent A first consequence of the evaluation formula is the
justification of the term ``virtual dimension'' for $d_\lambda$.

\Theorem virdim. Let $U$ be a multiplicity free space with ring of
functions $\cP(U)=\oplus_{\lambda\in\Lambda_+}\cP^\lambda$ and
associated $\rho$\_vector as in section~\cite{Fourier}. Then
$\|dim|\cP^\lambda=d_\lambda$ and $\|dim|U=2\ell(\rho)$.

\Proof: Let $\cD(U)=\oplus_{\lambda\in\Lambda_+}\cD_\lambda$ be the
decomposition of the space of constant coefficient differential
operators where $\cD_\lambda$ is simple with highest weight
$\lambda$. Fix $\lambda\in\Lambda_+$. If $D\in\cD_\lambda$ and
$f\in\cP^\lambda$ then $D(f)$ is a polynomial of degree zero, hence a
constant. This way we get a non\_degenerate pairing
$\cD_\lambda\times\cP^\lambda\pfeil\CC$. For any basis $f_i$ of
$\cP^\lambda$ let $D_i\in\cD_\lambda$ be the dual basis, i.e.,
$D_i(f_j)=\delta_{ij}$. Then $D:=\sum_if_iD_i$ is $G$\_invariant and
acts as identity on $\cP^\lambda$. By definition, the associated
polynomial $p_D$ is $p_\lambda$ (see \cite{Montreal}, or
\cite{Existence} \S7). We have $\tD=(-1)^{\ell(\lambda)}\sum_iD_if_i$,
hence $\tD(1)=(-1)^{\ell(\lambda)}\sum_i
D_i(f_i)=(-1)^{\ell(\lambda)}\|dim|\cP^\lambda$. On the other hand,
$p_\tD(v)=p_D(-v)=p_\lambda(-v)$ by \cite{Fourier-transform}. Thus
$\tD(1)=p_\tD(\rho)=p_\lambda(-\rho)=(-1)^{\ell(\lambda)}d_\lambda$
which shows $\|dim|\cP^\lambda=d_\lambda$.

The second formula is proved similarly. Here we choose a basis $x_i$
of $U^\vee\subseteq\cP(U)$. Let $\partial_i\in U\subseteq\cD$ be its
dual basis. Because $D=\sum_ix_i\partial_i$ is the Euler vector field
we have $p_D(z)=\ell(z-\rho)$. As above we get
$-\|dim|U=\tD(1)=p_D(-\rho)=-2\ell(\rho)$.\qed

\Remark: In the context of Hermitian symmetric
spaces the dimension formula was proved by Upmeier \cite{Up}.

\beginsection Interpolation. The interpolation formula

In this section we state a formula which allows to expand an arbitrary
$W$\_invariant polynomial in terms of the basis $p_\mu$. For this we
need another immediate consequence of the transposition formula
\cite{DualFormula}:

\Theorem Involution. Let $\rho\in V_0$ be dominant. Then the matrix
$$\eqno{E41}
\left((-1)^{\ell(\mu)}
p_\mu(\rho+\lambda)\right)_{\mu,\lambda\in\Lambda_+}
$$
is an involutory.

\Proof: By \cite{DualFormula}, the matrix expresses the involution
$h(z)\mapsto h(-z)$ of $\cP^W$ in the $q_\mu$\_basis.\qed

\noindent Let $\cC(\rho+\Lambda_+)$ be the set of $\CC$\_valued
functions on $\rho+\Lambda_+$. For $h\in\cC(\rho+\Lambda_+)$ we define
its {\it transform} $\widehat h\in\cC(\rho+\Lambda_+)$ by
$$\eqno{}
\widehat h(\rho+\mu):=\sum_{\tau\in\Lambda_+}(-1)^{\ell(\tau)}
p_\tau(\rho+\mu)h(\rho+\tau).
$$
The sum is finite since all summands with $\ell(\tau)>\ell(\mu)$ are
zero. We consider two subspaces of $\cC(\rho+\Lambda_+)$. First, let
$\cC_0(\rho+\Lambda_+)$ be the set of functions with finite
support. Secondly, we consider, via restriction, $\cP^W$ as subspace
of $\cC(\rho+\Lambda_+)$.

\Theorem Interpol. Let $\rho\in V_0$ be dominant. Then transformation
$h\mapsto\widehat h$ has the following properties:
\item{i)}$\widehat{\widehat h}=h$.
\item{ii)}$h\in\cP^W\Leftrightarrow\widehat h\in\cC_0(\rho+\Lambda_+)$.
\item{iii)}Interpolation formula:
$$\eqno{E11}
h(z)=\sum_{\mu\in\Lambda_+}(-1)^{\ell(\mu)}\widehat h(\rho+\mu)
p_\mu(z)
\quad\hbox{for all}\ h\in\cP^W.
$$\Par

\Proof: Let $a_{\mu\lambda}:=(-1)^{\ell(\mu)}p_\mu(\rho+\lambda)$. Then
$\widehat h(\rho+\mu)=\sum_\tau a_{\tau\mu}h(\rho+\tau)$ and therefore
$$\eqno{}
\widehat{\widehat h}(\rho+\lambda)=\sum_\mu a_{\mu\lambda}\widehat h(\rho+\mu)=
\sum_{\mu,\tau}a_{\mu\lambda}a_{\tau\mu}h(\rho+\tau)=
\sum_\tau\left[\sum_\mu a_{\tau\mu}a_{\mu\lambda}\right]h(\rho+\tau).
$$
By \cite{Involution}, the sum in brackets equals
$\delta_{\tau\lambda}$ which implies {\it i)}. 

Let $\chi_{\rho+\nu}\in\cC_0(\rho+\Lambda_+)$ be the characteristic
function of $\{\rho+\nu\}$. Then
$$\eqno{E13}
\widehat\chi_{\rho+\nu}=(-1)^{\ell(\nu)}p_\nu.
$$
Hence, $h\mapsto\widehat h$ maps a basis of $\cC_0(\rho+\Lambda_+)$ to
a basis of $\cP^W$ which proves {\it ii)}.

Finally, {\it i)} implies
$$\eqno{E10}
h(\rho+\lambda)=\sum_{\mu\in\Lambda_+}
(-1)^{\ell(\mu)}p_\mu(\rho+\lambda)\widehat h(\rho+\mu).
$$
By {\it ii)}, $\widehat h$ is a function with finite
support. Therefore, the sum \cite{E10} is over a finite set of $\mu$'s
which is independent of $\lambda$. This implies \cite{E11} since
$\rho+\Lambda_+$ is Zariski dense in $V$.\qed

The operator $L$ acts naturally on $\cC(\rho+\Gamma)$,
$$\eqno{}
(Lh)(\rho+\lambda)=\sum_\eta f_\eta(\rho+\lambda)h(\rho+\lambda-\eta),
$$
provided the coefficients $f_\eta(\rho+\lambda)$ are defined, i.e.,
$\rho$ is non\_integral. Then it follows from the first cut\_off
property of $f_\eta$, \cite{Cut-Off}{\it a)}, that the quotient
$\cC(\rho+\Lambda_+)$ is $L$\_stable. Let $\cA$ be the algebra
generated by $\cP^W$ and $L$ in $\|End|\cC(\rho+\Lambda_+)$. It
follows that $\cC(\rho+\Lambda_+)$ is an $\cA$\_module. Moreover,
$\cC_0(\rho+\Lambda_+)$ and $\cP^W$ are $\cA$\_submodules.  For every
$X\in\|End|_\CC\cC(\rho+\lambda_+)$ we define $\widehat X$ by
$\widehat X(h):=\widehat{X(\widehat h)}$.

\Theorem InvAut. Assume $\rho$ is non\_integral. Then $X\mapsto
\widehat X$ induces an involutory automorphism of $\cA$. More
precisely, we have $\widehat m_h=D_h$ and $\widehat L=-L$. Here, $m_h$
is the operator\footnote{\rm Since $\widehat m_h\ne m_{\widehat h}$ we
are forced to use this notation.} ``multiplication by $h$'' .

\Proof: The equality $\widehat m_h=D_h$ is equivalent to
$$\eqno{E12}
\widehat{h\widehat p}(\rho+\lambda)=D_h(p)(\rho+\lambda)
\quad\hbox{for all }p\in\cC(\rho+\Lambda),\lambda\in\Lambda_+.
$$
Now we fix $\lambda$. Then both sides of \cite{E12} depend only on the
values of $p$ in finitely many points, more precisely, in points
$\rho+\mu$ with $\ell(\mu)\le\ell(\lambda)$. Since there is a
$W$\_invariant polynomial which has the same values at these points we
may assume $p\in\cP^W$. By linearity, we may assume $p=p_\nu$. Then,
by \cite{E13},
$$\eqno{}
h\widehat p_\nu=(-1)^{\ell(\nu)}h\chi_{\rho+\nu}=
(-1)^{\ell(\nu)}h(\rho+\nu)\chi_{\rho+\nu}=h(\rho+\nu)\widehat p_\nu
$$
and therefore
$$\eqno{}
\widehat m_h(p_\nu)=h(\rho+\nu)\widehat{\widehat p_\nu}=
h(\rho+\nu)p_\nu=D_h(p_\nu).
$$
This proves $\widehat m_h=D_h$. But then $\widehat
L=(\ell-D_\ell)^\wedge=D_\ell-\ell=-L$. This shows in particular
that $X\mapsto \widehat X$ maps $\cA$ into itself.\qed

\Remark: The non\_integrality of $\rho$ is needed to make sense of the
action of $\cA$ on $\cC(\rho+\Lambda_+)$. As already mentioned, the
element $\rho$ attached to a multiplicity free representation is never
non\_integral. It will be a consequence of \cite{effect} that
$X\mapsto\widehat X$ is, in fact, defined for every $\rho\in V_0$.
\medskip

\beginsection Product. The scalar product

Assume $\rho$ is strongly dominant. The symmetry property
\cite{Symmetrisch} indicates the presence of a scalar product on
$\cP^W$. In fact, we define a non\_degenerate scalar product on $\cP^W$ by
$$\eqno{E16}
\<p_\lambda,p_\mu\>=d_\lambda\delta_{\lambda\mu}
\quad\hbox{for all }\lambda,\mu\in\Lambda_+.
$$
Thus $\<p_\lambda,q_\mu\>=(-1)^{\ell(\lambda)}\delta_{\lambda\mu}$ and
$\<q_\lambda,q_\mu\>=d_\lambda^{-1}\delta_{\lambda\mu}$. For any
function $h(z)$ let $h^-(z):=h(-z)$. Then we have

\Theorem. For all $\lambda\in\Lambda_+$ and $h\in\cP^W$ holds
$$\eqno{E15}
\<q_\lambda,h\>=\widehat h(\rho+\lambda)\hbox{ and }
\<q_\lambda^-,h\>=h(\rho+\lambda).
$$

\Proof: From the interpolation formula \cite{E11} we obtain
$$\eqno{}
\<q_\lambda,h\>=\sum_\mu(-1)^{\ell(\mu)}\widehat
h(\rho+\mu)\<q_\lambda,p_\mu\>=\widehat h(\rho+\lambda).
$$
Moreover, from \cite{DualFormula}, \cite{E15}, and \cite{E11} we get
$$\eqno{}
\<q_\lambda^-,h\>=\sum_\mu(-1)^{\ell(\mu)}p_\mu(\rho+\lambda)\<q_\mu,h\>=
\sum_\mu(-1)^{\ell(\mu)}\widehat
h(\rho+\mu)p_\mu(\rho+\lambda)=h(\rho+\lambda).
$$
\qed

\Remark: In particular, we have
$\<q_\lambda^-,q_\mu^-\>=q_\mu(-\rho-\lambda)$ which explains the
symmetry in $\lambda$ and $\mu$.
\medskip

\noindent
There is also a general expression for the scalar product:

\Theorem. For all $g,h\in\cC^W$:
$$\eqno{E42}
\<g,h\>=\sum_{\mu\in\Lambda_+}d_\mu\,\widehat g(\rho+\mu)\widehat h(\rho+\mu).
$$

\Proof: Just apply the interpolation formula \cite{E11} to $g$ and $h$.\qed

The algebra $\cA$ is not quite closed under taking adjoints for the
scalar product. Therefore, these will be studied in the next
section. Here, we use a slightly modified scalar product:
$$\eqno{E14}
\<g,h\>^-:=\<g^-,h^-\>.
$$
The adjoint of an operator $X$ with respect to the scalar product
\cite{E14} will be denoted by $X'$.

\Theorem Adjoint1. Let $\rho\in V_0$ be strongly dominant. Then for
every $X\in\cA$ the adjoint $X'$ exists and is again in $\cA$. More
precisely, $h'=D_{h^-}$, and $L'=L$. In particular, $X\mapsto X'$
induces an involutory antiautomorphism of $\cA$. Moreover, $(\widehat
X)'=(X')^\wedge$ for all $X\in\cA$.

\Proof: By \cite{E15} we have $\<q_\lambda,h\>^-=h(-\rho-\lambda)$ for all
$h\in\cP^W$. Then I claim
$$\eqno{}
\<D_{h^-}(f),g\>^-=\<f,hg\>^-.
$$
for all $f,g,h\in\cP^W$. Indeed, it suffices to prove this for
$f=q_\lambda$. Then
$$\eqno{}
\<D_{h^-}(q_\lambda),g\>^-=h(-\rho-\lambda)\<q_\lambda,g\>^-=
h(-\rho-\lambda)g(-\rho-\lambda)=\<q_\lambda,hg\>^-.
$$
Thus the adjoint operator of $h$ is $D_h$. Then we also have
$$
L'=(\ell-D_\ell)'=(-\ell^--D_\ell)'=-D_\ell-\ell^-=L.
$$
Finally, $(\widehat L)'=-L=(L')^\wedge$ and $(\widehat m_h)'=
D_h'=m_{h^-}=\widehat D_{h^-}=(m_h')^\wedge$ which shows the last
claim.\qed

\beginsection sl-two. The $PGL_2$\_action

For any operator $X\in\|End|_\CC(\cP^W)$ define the operator $X^-$ by
$X^-(g)=X(g^-)^-$. In particular, if $X=\sum_\tau a_\tau(z)T_\tau$ is
a difference operator then $X^-=\sum_\tau a_\tau(-z)T_{-\tau}$ is
again a difference operator. For multiplication operators we have
$m_h^-=m_{h^-}$. On the other side, $L^-$ is new. Therefore, let $\cB$
be the algebra generated by $\cP^W$, $L$, and $L^-$. It contains $\cA$
as a subalgebra. Moreover $X\mapsto X^-$ induces an involutive
automorphism of $\cB$. Observe that $\cA$ contains only operators
composed of shifts by $\tau\in\Lambda$ while in $\cB$ arbitrary shifts
$\tau\in\Gamma^\vee$ are possible.

For any $X\in\|End|_\CC(\cP^W)$ let $X^*$ be the adjoint operator (if
it exists) with respect to the scalar product $\<\cdot,\cdot\>$ defined
in \cite{E16}. Its relation to the adjoint $X'$ is
$X^*=X^{-\prime-}$. Indeed
$$\eqno{}
\<Xf,g\>=\<X^-f^-,g^-\>^-=\<f^-,X^{-\prime}g^-\>^-=\<f,X^{-\prime-}g\>.
$$

\Theorem adjoints. Let $\rho\in V_0$ be strongly dominant. Then for
every $X\in\cB$ the adjoint operator $X^*$ exists and is again in
$\cB$. More precisely, the following formulas hold (with $h\in\cP^W$):
$$\eqno{E17}
\vcenter{\halign{$#$\hfill&$#$\hfill\cr
h^*&=D_h^-=\|exp|(\|ad|L^-)(h^-),\cr
L^*&=E^--E=L-2\ell-L^-,\cr
(L^-)^*&=L^-,\cr
D_h^*&=D_h.\cr}}
$$
In particular, $X\mapsto X^*$ induces an involutive antiautomorphism
of $\cB$.

\Proof: Since $D_h$ has an orthogonal eigenbasis, $p_\lambda$, it is
self\_adjoint: $D_h^*=D_h$. By \cite{Adjoint1} we have
$h^*=((h^-)')^-=D_h^-$. Moreover,
$$\eqno{}
L^*=(\ell-D_\ell)^*=D_\ell^--D_\ell=
E^--E=(\ell-L)^--(\ell-L)=L-2\ell-L^-.
$$
Finally, $(L^-)^*=L^{\prime-}=L^-$.\qed

\Remark: Of course, $\cB$ is still preserved under the other adjoint
$X\mapsto X'$ with $(L^-)'=(L^-)^{-*-}=L^{*-}=E-E^-=-L^-$.
\medskip
Recall that three elements $(e,h,f)$ of a (Lie) algebra are called an
$sl_2$\_triple if the relations $[h,e]=2e$, $[h,f]=-2f$, and $[e,f]=h$
hold.

\Theorem triple. Both $(L,2\ell,L^-)$ and $(-L,2E,L^*)$ are
$sl_2$\_triples.

\Proof: For every $\eta\in\Gamma^\vee$ holds
$[\ell,T_\eta]=\ell(\eta)T_\eta$. Hence, by definition of $L$, we have
$[2\ell,L]=2L$. We also get $[2E,L]=[2\ell-2L,L]=2L$. The equation
$[2\ell,L^-]=-2L^-$ follows by applying $X\mapsto X^-$ to both sides
of $[2\ell,L]=2L$. Moreover, if we apply $X\mapsto X^*$ to $[2E,L]=2L$
we get, according to \cite{E17}, $[2E,L^*]=-[2E,L]^*=-2L^*$. Moreover,
$$\eqno{}
\eqalign{
[L,L^-]&=[L,L-2\ell-L^*]= 2L-[\ell-E,L^*]=\cr
&=2L-[\ell,L-2\ell-L^-]-L^*= 2L-L-L^--(L-2\ell-L^-)=2\ell.\cr}
$$
Finally,
$[-L,L^*]=[-L,L-2\ell-L^-]=-2L+2\ell=2E$.\qed

\noindent Of course, the two triples span the same three dimensional
subspace $\fs$ inside $\cB$ which we identify with the Lie algebra
$sl_2(\CC)$ by using the second triple:
$$\eqno{}
-L\mapsto\pmatrix{0&1\cr0&0\cr},\quad
2E\mapsto\pmatrix{1&0\cr0&-1\cr},\quad
L^*\mapsto\pmatrix{0&0\cr1&0\cr}.
$$
Then we also have
$$\eqno{}
2\ell\mapsto\pmatrix{1&-2\cr0&-1\cr},\quad
L^-\mapsto\pmatrix{-1&1\cr-1&1\cr}.
$$

\def\bmatrix#1{\left[\matrix{#1}\right]}
Now we would like to integrate the inner $\fs$\_action on $\cB$. For
this, let $S:=\|Aut|\fs$. Its Lie algebra is $\fs$. Moreover, if we
identify $\fs$ with $sl_2(\CC)$ as above then $S$ gets identified with
$PGL_2(\CC)$. Its elements are invertible $2\times2$\_matrices modulo
scalar multiplication which we write in square brackets. Of particular
interest is the involution
$$\eqno{}
\sigma:=\bmatrix{1&-1\cr0&-1\cr}\in S
$$
which maps the two $sl_2$\_triples into each other: 
$$\eqno{}
(L,2\ell,L^-)=\sigma(-L,2E,L^*).
$$

\Theorem integration. The adjoint action of $\fs$ on $\cB$ can be
integrated to an algebraic $S$\_action.

\Proof: First, we show that $\|ad|\fs$ acts locally finitely on
$\cB$. By Poincar\'e\_Birkhoff\_Witt it suffices to show that
for $L$, $2\ell$, and $L^-$, separately.

We claim that the elements $L$ and $L^-$ act locally nilpotently. It
suffices to show this on the generators $h\in\cP^W$, $L$, and
$L^-$. For $L$, the assertion follows from \cite{nilp} (for $h$) and
\cite{triple} (for $L^-$). For $L^-$ we apply the automorphism
$X\mapsto X^-$.

The action of $\|ad|2\ell$ on difference operators is clearly
diagonalizable. This shows already that $\|ad|\fs$ integrates to an
$SL_2(\CC)$\_action. The possible eigenvalues of $\|ad|2\ell$ are
$2\ell(\tau)$, $\tau\in\Gamma^\vee$. Since these are all even, the
action of $SL_2(\CC)$ descends to an action of $PGL_2(\CC)=S$.\qed

\noindent Now we compute the effect of some particular elements of $S$
on $\cB$.

\Proposition effect. The effect of $\sigma$ on the generators of $\cB$ are
$$\eqno{}
\sigma(L)=-L,\quad\sigma(h)=D_h,\quad\sigma(L^-)=L^*.
$$

\Proof: We already know $(L,2\ell,L^-)=\sigma(-L,2E,L^*)$. Thus it
remains to calculate $\sigma(h)$. To this end, write
$\sigma=\alpha\beta$ where
$$\eqno{}
\alpha=\bmatrix{1&-1\cr0&1\cr},\quad\beta=\bmatrix{1&-2\cr0&-1\cr}=
\sigma\bmatrix{1&0\cr0&-1\cr}\sigma^{-1}.
$$
The matrix $\beta$ lies in the Cartan subgroup whose Lie algebra is
$\CC\ell$. Therefore, it fixes every element of $\cB$ which commutes
with $\ell$. This implies $\beta(h)=h$. The matrix $\alpha$ acts by
$\|exp|(\|ad|L)$ on $\cB$. Hence it sends, by definition, $h$ to
$D_h$. We conclude $\sigma(h)=D_h$.\qed

\noindent Next we investigate the effect of $\fs$ on the $\cB$\_module
$\cP^W$.

\Theorem LL*. Let $\rho\in V_0$ be strongly dominant. Then for all
$\lambda\in\Lambda_+$ and $d\in\NN$ holds
$$\eqno{L-Formula}
{1\over d!}L^d(p_\lambda)=
\sum_{\mu\in\Lambda_+\atop\ell(\mu)=\ell(\lambda)+d}
p_\lambda(\rho+\mu)p_\mu,
$$
$$\eqno{L*-Formula}
{1\over d!}(-L^*)^d(q_\lambda)=
\sum_{\mu\in\Lambda_+\atop\ell(\mu)=\ell(\lambda)-d}
p_\mu(\rho+\lambda)q_\mu.
$$

\Proof: By \cite{triple} we have $[E,L]=L$, hence $[E,L^d]=dL^d$. For
every $\lambda\in\Lambda_+$ follows that $L^d(p_\lambda)$ is a linear
combination of those $p_\mu$ with $\ell(\mu)=\ell(\lambda)+d$. On the
other hand, we have
$L^d(p_\lambda)=(\ell-E)^d(p_\lambda)=\ell^dp_\lambda$ plus lower
order terms. Then \cite{L-Formula} follows from \cite{Pieriell}.

Using the fact that the dual basis of the $p_\lambda$ are the
$(-1)^{\ell(\lambda)}q_\lambda$ we get from \cite{L-Formula}
$$\eqno{L*-Formula2}
{1\over d!}(L^*)^d((-1)^{\ell(\lambda)}q_\lambda)=
\sum_{\mu\in\Lambda_+\atop\ell(\mu)=\ell(\lambda)-d}
p_\mu(\rho+\lambda)(-1)^{\ell(\mu)}q_\mu
$$
which is equivalent to \cite{L*-Formula}.\qed

\noindent Formulas \cite{L-Formula} and \cite{L*-Formula} can be
expressed more conveniently as generating series:
$$\eqno{E20}
\|exp|(tL)p_\lambda=\sum_{\mu\in\Lambda_+}
t^{\ell(\mu)-\ell(\lambda)}p_\lambda(\rho+\mu)p_\mu
$$
and
$$\eqno{E21}
\|exp|(-tL^*)q_\lambda=\sum_{\mu\in\Lambda_+}
t^{\ell(\lambda)-\ell(\mu)}p_\mu(\rho+\lambda)q_\mu.
$$
There is a big difference between this two formulas in that the
latter, \cite{E21}, is a finite sum. This means that \cite{E21}
defines an algebraic action of
$$\eqno{}
\bmatrix{1&0\cr-t&1\cr}
$$
on $\cP^W$. There is also an action of the diagonal matrices on $\cP$,
defined by
$$\eqno{E22}
\bmatrix{a&0\cr0&b\cr}:
q_\lambda\mapsto\left({a\over b}\right)^{\ell(\lambda)}q_\lambda.
$$
Then \cite{E21} and \cite{E22} combine to an action of $B$, the
subgroup of lower triangular matrices of $S=PGL_2(\CC)$. This action
is compatible with that on $\cB$:
$$\eqno{}
{}^b(Xh)={}^bX(\,{}^bh)\quad\hbox{for all }b\in B,X\in\cB,h\in\cP^W.
$$

\Remark: The action of $B$ is on $\cP^W$ is not quite the one which
one would obtain by exponentiating the action of
$\|Lie|B\subset\fs\subset\cB$ on $\cP^B$. The reason is that
$q_\lambda$ is an eigenvector of $E$ with eigenvalue
$\ell(\lambda)+\ell(\rho)$ and not just $\ell(\lambda)$. Therefore,
unless $\ell(\rho)$ is an integer, the exponentiated $B$\_action is
not algebraic. In the geometric case, i.e., when $\rho$ comes from a
multiplicity free action on a vector space $U$, we have that
$\ell(\rho)={1\over2}\|dim|U$ (\cite{virdim}) is in ${1\over2}\ZZ$. In
that case, one can integrate the $\|Lie|B$\_action to an algebraic
action of the lower triangular matrices in $SL_2(\CC)$.  \medskip

\noindent Now we can locate the automorphism $X\mapsto X^-$ in
$S$:

\Theorem minusaction. The matrix
$$\eqno{}
\gamma:=\bmatrix{1&0\cr1&-1\cr}\in B
$$
acts as $h\mapsto h^-$ on $\cP^W$ and as $X\mapsto X^-$ on $\cB$.

\Proof:
We write $\gamma=\alpha\beta$ with
$$\eqno{}
\alpha=\bmatrix{1&0\cr 1&1\cr},\quad\beta=\bmatrix{1&0\cr0&-1\cr}.
$$
Then $\beta(q_\lambda)=(-1)^{\ell}q_\lambda$ (by~\cite{E22}) and
$\alpha(q_\lambda)=
\sum_\mu(-1)^{\ell(\lambda)-\ell(\mu)}p_\mu(\rho+\lambda)q_\mu$
(by~\cite{E21}).  The transposition formula~\cite{DualFormula} implies
$\gamma(q_\lambda)=q_\lambda^-$. We conclude $\gamma(h)=h^-$ by
linearity. Finally,
$\gamma(X)(h)=\gamma(X(\gamma(h))=X(H^-)^-=X^-(h)$.\qed

\Remark: One consequence of \cite{minusaction} is the formula 
$$\eqno{}
\|exp|(L^*)(p_\lambda)=(-1)^{\ell(\lambda)}p_\lambda^-.
$$
It has the advantage that it works for $\rho$ which are just dominant.
\medskip
Now we come back to the automorphism $X\pfeil\widehat X$ of
section~\cite{Interpolation}. Comparing \cite{InvAut} with
\cite{effect} we see that $\sigma$ induces on $\cA$ exactly
$X\pfeil\widehat X$. Now we extend this to $\cB$:

\Theorem. Let $\rho$ be non\_integral. Then $\cC(\rho+\Lambda_+)$ is
naturally a $\cB$\_module. Moreover the relation $(L^-)^\wedge=L^*$
holds. In particular, we have $\widehat X=\sigma(X)$ for all $X\in\cB$.

\Proof: By definition, we have
$$\eqno{}
L^-(h)(\rho+\mu)=\sum_{\eta\in\Lambda_1}f_\eta(-\rho-\mu)h(\rho+\mu+\eta).
$$
Thus, it follows from \cite{Cut-Off}{\it b)} that $L^-$ and therefore
$\cB$ acts on $\cC(\rho+\Lambda_+)$.

For every fixed $\lambda$ the values $(L^-)^\wedge(h)(\rho+\lambda)$
and $L^*(h)(\rho+\lambda)$ depend on only finitely many values of
$h$ which we may interpolate by a linear combination of
$p_\lambda$'s. This implies, that it suffices to prove
$(L^-)^\wedge(h)=L^*(h)$ for $h=p_\lambda$. We have
$$\eqno{}
\eqalign{(L^-)^\wedge(p_\lambda)&=(-1)^{\ell(\lambda)}(L^-)^\wedge(\widehat\chi_{\rho+\lambda})=(-1)^{\ell(\lambda)}(L^-\chi_{\rho+\lambda})^\wedge=\cr
&=(-1)^{\ell(\lambda)}
\sum_\eta f_\eta(-\rho-\lambda+\eta)\widehat\chi_{\rho+\lambda-\eta}
=-\sum_\eta f_\eta(-\rho-\lambda+\eta)p_{\lambda-\eta}\cr}
$$
Since $D_\ell=\ell-L$ we have $f_\eta(z)=-b_\eta^\ell(z)$. Therefore,
if we compare \cite{P-Formula-1} (with $h=\ell$) and \cite{Pieriell}
(with $k=1$) we get
$$\eqno{}
-f_\eta(-\rho-\lambda+\eta)=
{d_\lambda\over d_{\lambda-\eta}}p_{\lambda-\eta}(\rho+\lambda)
$$
Thus, using \cite{L*-Formula} (with $d=1$) we get
$$\eqno{}
\eqalign{
(L^-)^\wedge(p_\lambda)&=(-1)^{\ell(\lambda)-1}d_\lambda\sum_\eta
p_{\lambda-\eta}(\rho+\lambda)q_{\lambda-\eta}=\cr
&=-(-1)^{\ell(\lambda)}d_\lambda (-L^*)(q_\lambda)=L^*(p_\lambda).\cr}
$$\qed

\beginsection differential. The differential limit

In this section we consider the effect of our difference operators on
the highest degree component of a polynomial. Let $\cP_{\le
d}:=\{h\in\cP\mid \|deg|h\le d\}$, $\cPq_d:=\cP_{\le d}/\cP_{\le d-1}$
and $\cPq:=\oplus_d\cPq_d$, the associated graded algebra. Observe
that $\cPq\cong\cP$ (even equivariantly) since $\cP$ is a polynomial ring.

Now we introduce the degree of an operator $X\in\cB$ as
$$\eqno{}
\|deg|X:=\|max|\,\{\|deg|X(h)-\|deg|h\mid h\in\cP^W\}.
$$
It is clear that the degree of any difference operator is finite.  Let
$\cB_{\le d}:=\{X\in\cB\mid\|deg|X\le d\}$. This defines a filtration
of $\cB$, i.e., $\cB_d$ is a subspace of $\cB$ with $\cB=\cup_d\cB_d$
and $\cB_d\cB_e\subseteq\cB_{d+e}$. Let $\cBq_d:=\cB_{\le d}/\cB_{\le
d-1}$ and $\cBq:=\oplus_d\cBq_d$, the associated graded algebra. The
point is now that more or less by construction, $\cPq^W$ is a faithful
$\cBq$\_module. We call it the {\it differential limit\/} since:

\Proposition diffop. Every $\Xq\in\cBq$ acts as a differential
operator on $\cPq^W$.

\Proof: We may assume that $\Xq\in\cBq_d$ is non\_zero and that it is
represented by a difference operator $X\in\cB_{\le d}$. Choose linear
coordinates $z_1,\ldots,z_n\in V^\vee$. By Taylor's theorem, the
translation operator $T_\eta$ can be written as differential operator
of infinite order:
$$\eqno{}
T_\eta=\|exp|\left(-\sum_i z_i(\eta){\partial\over\partial z_i}\right).
$$
Therefore, we can also expand $X$ into an infinite order differential
operator with coefficients of bounded degree.

Now let $h\in\cP^W$ be a
polynomial of degree $e$ with highest degree component
$\hq$. For an indeterminate $t$ let
$h_t(z):=h(t^{-1}z)$. Then
$$\eqno{}
h_t=\hq t^{-e}+\ldots
$$
where ``$\ldots$'' means ``terms of higher order in $t$''.

Correspondingly, we define $X_t$ by $X_t(h):=X(h_{t^{-1}})_t$. This
amounts to replacing all variables $z_i$ by $t^{-1}z_i$ and all
partial derivatives ${\partial\over\partial z_i}$ by
$t{\partial\over\partial z_i}$. In particular, we have
$X_t(h_t)=X(h)_t$. Now we develop $X_t$ into a Laurent series in
$t$. This is possible since the coefficients of $X$ have bounded
degree. Thus there is $N\in\ZZ$ with
$$\eqno{}
X_t=\XS t^{-N}+\ldots
$$
where $\XS$ is a non\_zero differential operator. Hence
$$\eqno{}
X(h)_t=X_t(h_t)=\XS(\hq)t^{-e-N}+\ldots
$$
This shows that $\|deg|X(h)\le\|deg|h+N$ with equality for most
$h$. Therefore, $N=d$ and $\Xq(\hq)=\XS(\hq)$. Thus $\Xq=\XS$ is a
differential operator.\qed

For the reminder of this section we assume that $\rho$ is dominant. We
show that the pair $(\cB,\cP^W)$ is isomorphic to $(\cBq,\cPq^W)$. For
this we use the action of the difference Euler operator $E$. Its
action on $\cP^W$ is diagonalizable with eigenvalues of the form
$d+\ell(\lambda)$, $d\in\NN$. Therefore, let $\cP_d^W:=\{h\in\cP^W\mid
E(h)=(d+\ell(\rho))h\}$. Then $\cP^W=\oplus_d\cP^W_d$. A basis of
$\cP^W_d$ is formed by all $p_\lambda$ with $\ell(\lambda)=d$. Thus
\cite{Main-Prop-E} implies that
$$\eqno{}
\cP^W_{\le d}=\oplus_{i\le d}\cP^W_i.
$$
In particular, the projection $\cP^W_d\pfeil\cPq^W_d$ is an isomorphism. This
way, we get an isomorphism (of vector spaces)
$$\eqno{}
\psi:\cP^W=\oplus_d\cP^W_d\Pf\sim\oplus_d\cPq_d^W=\cPq^W.
$$

Now we do the same thing with $\cB$. We know from the last section
that the action of $\|ad|E$ on $\cB$ is diagonalizable with integral
eigenvalues. Therefore, let $\cB_d:=\{X\in\cB\mid [E,X]=dX\}$. Then
$\cB=\oplus_d\cB_d$ is a grading of $\cB$.

\Lemma. Let $\rho$ be dominant. Then $\cB_{\le d}=\oplus_{i\le d}\cB_i$.

\Proof: Let $X\in\cB_d$ and $h\in\cP_e$. Then
$$\eqno{}
EX(h)=[E,X](h)+XE(h)=dX(h)+eX(h)=(d+e)X(h)
$$
implies $\cB_d\cP_e\subseteq\cP_{d+e}$. In particular, we have
$\cB_d\cP_{\le e}\subseteq\cP_{\le d+e}$ which shows $\cB_{\le
d}\supseteq\oplus_{i\le d}\cB_i$.

Conversely, let $X\in\cB_{\le d}$ and $X=\sum X_n$ with $X_n\in\cB_n$
and $N=\|max|\{n\mid X_n\ne0\}$. Choose $h\in\cP_e$ with
$X_N(p)\ne0$. Since $X_i(p)\subseteq\cP_{i+e}$ is either zero or has
precisely the degree $i+e$ we conclude $\|deg|X(p)=N+e$. From the
assumption $\|deg|X(p)\le d+e$ follows $N\le d$. This proves
$X\in\oplus_{i\le d}\cB_i$.\qed

\noindent An immediate consequence of the lemma is
$\cB_d\pf\sim\cBq_d$ which gives rise to a map
$$\eqno{E23}
\Psi:\cB=\oplus_d\cB_d\Pf\sim\oplus_d\cBq_d=\cBq.
$$

\Theorem difflim. Let $\rho$ be dominant. Then the map $\Psi$ in
\cite{E23} is an isomorphism of algebras. Moreover, under this
isomorphism the $\cB$\_module $\cP^W$ corresponds to the
$\cBq$\_module $\cPq^W$. More precisely,
$$\eqno{}
\psi(Xh)=\Psi(X)\psi(h)\quad\hbox{for all }X\in\cB,h\in\cP^W.
$$

\Proof: The relation $\Psi(XY)=\Psi(X)\Psi(Y)$ has to be proven only
for $X\in\cB_d$, $Y\in\cB_e$. But then it follows from
$\cB_d\cB_e\subseteq\cB_{d+e}$. Similarly, for \cite{E23} we may
assume $X\in\cB_d$ and $h\in\cP_e$. Then it follows from
$\cB_d\cP_e\subseteq\cP_{d+e}$.\qed

\noindent In view of this theorem it is probably more adequate to call
$\cBq$ the differential ``picture'' as opposed the differential
``limit'' of $\cB$. It shows that the difference operators are just
represented differently namely by differential operators.

Next, we study the maps $\psi$ and $\Psi$ more closely. Given
$h\in\cP^W$, there are two ways to produce an element of $\cPq^W$:
first $\hq$, its top homogeneous component, and then $\psi(h)$. We
have $\psi(h)=\hq$ precisely if $h$ is an $E$\_eigenvector. Therefore,
consider $\pq_\lambda$, the top homogeneous component of
$p_\lambda$. These polynomials are also of high representation
theoretic interest (see, e.g., \cite{Montreal}. In the classical case
they are the Jack polynomials.) They form a basis of
$\cPq^W$. Since $p_\lambda$ is an $E$\_eigenvector we could {\it define}
$\psi$ by the property $\psi(p_\lambda)=\pq_\lambda$.

The same thing works for $\cB$: every $X\in\cB$ gives rise to two
elements in $\cBq$ namely its top homogeneous component $\Xq$ and
$\Psi(X)$. Moreover, $\Psi(X)=\Xq$ if and only if $X$ is an
$\|ad|E$\_eigenvector. This holds in particular for $\cB_0$, the
commutant of $E$. Hence we have $\Psi(D_h)=\Dq_h$ where $\Dq_h$ are
certain differential operators. In the classical case, they are the
Sekiguchi\_Debiard operators, \cite{Sek}, \cite{Deb}. They are
simultaneously diagonalized by the $\pq_\lambda$:
$$\eqno{}
\Dq_h(\pq_\lambda)=
h(\rho+\lambda)\pq_\lambda\quad\hbox{for all }h\in\cP^W.
$$

Next, we compute the image of the $sl_2$\_triple $(-L,2E,L^*)$.

\Proposition. We have $\Psi(L)=m_{\overline\ell}$ (multiplication by
$\overline\ell\in\cPq^W$) and $\Psi(E)=\Eq=\xi+\ell(\rho)$ where $\xi$
is the Euler vector field. The differential operator
$\Lq^*:=\Psi(L^*)$ is of order 2 and of degree $-1$.

\Proof: We have $L\in\cB_1$, hence $\Psi(L)=\Lq$. From
$\|deg|(m_\ell-L)=\|deg|E<1$ it follows
$\Lq=\overline{m_\ell}=m_{\overline\ell}$. Since $\Eq$ acts on
$\cq^W_d$ by multiplication with $d+\ell(\rho)$ we have
$\Psi(E)=\Eq=\xi+\ell(\rho)$. Since $L^*\in\cB_{-1}$, the degree of
$\Lq^*$ is $-1$. Expand $L^*_t$ as a Laurent series in $t$ as in the
proof of \cite{diffop}. Since the coefficients of $L^*=L-2\ell-L^-$
are rational functions of degree 1 we have
$$\eqno{}
L^*_t=X_0t^{-1}+X_1+X_2t+\ldots
$$
where $X_i$ is homogeneous of degree $1-i$. Thus $X_0=X_1=0$ and
$X_2=\Lq*$. By construction, $X_i$ is a differential operator of
order $i$. Therefore, the order of $\Lq^*$ is $2$.\qed

Now we compare the multiplication operators in $\cB$ and $\cBq$.

\Theorem. a) Let $h\in\cP^W$. Then
$$\eqno{}
\Psi(m_h)=\|exp|(-\|ad|\overline\ell)(\Dq_h).
$$
b) Conversely, let $\hq\in\overline\cP^W_d$ and choose a lift $h\in\cP^W_{\le
d}$. Then
$$\eqno{E25}
\Psi^{-1}(m_\hq)={1\over d!}(-\|ad|L)^d(m_h).
$$

\Proof: a) We have
$$\eqno{E24}
m_h=\|exp|(-\|ad|L)(D_h)=\sum_i{1\over i!}(-\|ad|L)^i(D_h).
$$
Each summand is $\|ad|E$\_homogeneous. Therefore,
$$\eqno{}
\Psi(m_h)=\sum_i{1\over d!}\overline{(-\|ad|L)^d(D_h)}=
\sum_i{1\over d!}(-\|ad|\Lq)^d(\Dq_h)=\|exp|(-\|ad|\overline\ell)(\Dq_h).
$$
b) Let $R$ denote the right hand side of \cite{E25}. The sum in
\cite{E24} terminates at $i=d$. Moreover, the $i$\_th summand is
$\|ad|E$\_homogeneous of degree $i$. This implies
$\Psi(R)=\Rq=\overline{m_h}=m_\hq$.\qed

\noindent Thus we obtained besides the $m_h$ and the $D_h$ yet another
commutative subalgebra of $\cB$ formed by the $\Psi^{-1}(m_\hq)$.

Finally, we discuss the geometric situation: let $U$ be a multiplicity
free space as in section~\cite{Fourier}. Since $\cP^W$ can be
identified with the algebra of $G$\_invariant differential operators
on $U$ we can use the symbol map to identify $\cPq^W$ with the algebra
of $G$\_invariant functions on the cotangent bundle, i.e., with
$\cP(U\oplus U^\vee)^G$. On the other hand we can think of $\cPq^W$ as
$W$\_invariant functions on $V$, i.e., of functions on $V/W$.

Now consider $\cP\cD(U\oplus U^\vee)^G$, the algebra of $G$\_invariant
differential operators on $U\oplus U^\vee$. These act on
$G$\_invariants and therefore we get a map
$$\eqno{}
\Phi:\cP\cD(U\oplus U^\vee)^G\pfeil\cP\cD(V/W).
$$
(This is an analogue of the Harish Chandra homomorphism.)  Observe
that $\cBq\subseteq\cP\cD(V/W)$.

\Theorem. The algebra $\cBq$ is in the image of $\Phi$.

\Proof: The algebra $\cB$ is generated by $L$, $\{m_h\mid
h\in\cP^W\}$, and $L^*$. Because of \cite{E24} we can replace $m_h$ by
$D_h$. Applying $\Psi$, we see that $\cBq$ is generated by
$m_{\overline\ell}$, $\{\Dq_h\mid h\in\cP^W\}$, and $\Lq^*$. We show
that these generators lie in the image of $\Phi$.

Choose coordinates $(x_1,\ldots,x_n,y_1,\ldots,y_n)$ of $U\oplus
U^\vee$ such that the natural pairing between $U$ and $U^\vee$ is
given by $q:=\sum_ix_iy_i$. Then $q$ is the symbol of the Euler vector
field and therefore $\Phi(q)=\overline\ell$.

We have $\cP\cD(U)^G\into\cP\cD(U\oplus U^\vee)^G$ by letting
operators act on the first factor. Thus we have a map
$\cP^W\pfeil\cP\cD(V/W)$ whose image are the differential operators
$\Dq_h$ (see \cite{Montreal} Thm.~4.11).

Finally, let $\Delta:=-\sum_i{\partial^2\over\partial x_i\partial
y_i}$ be the Laplace operator. Then it follows from
\cite{BenRat} (1.8)\footnote{The $p_\alpha$ in that paper is our
$\qq_\alpha$.} that $\Phi(\Delta)$ acts on the $\qq_\lambda$ exactly
as $\Lq^*$.\qed

\noindent{\bf Question:} Is the image of $\Phi$ exactly $\cBq$?

\beginsection Binomial. The binomial formula

In this section we investigate another limiting case, namely, we are
looking at the infinitesimal neighborhood of a point
$\delta\in\Sigma^\vee\cap V^W$. We are going to prove a binomial type
formula for $\pq_\lambda(z+\delta)$.

The set $\Sigma^\vee\cap V^W$ has usually just one element but there
are cases\footnote{Cases {\bf III} ($n$ odd), {\bf IVa}, and {\bf IVc}
of \cite{Existence}~\S8} where it is empty and there is one
case\footnote{Case {\bf V} of \cite{Existence}~\S8} where it consists
of two points.  For the classical or semiclassical case see the
example below.

Let $\omega_\delta\in\Sigma$ be the dual element for $\delta$, i.e.,
$$\eqno{}
\omega(\delta)=\cases{1&if $\omega=\omega_\delta$,\cr
0&if $\omega\in\Sigma$ and $\omega\ne\omega_\delta$.\cr}
$$
Since $\delta$ is $W$\_invariant we have $\omega(\delta)\in\{0,1\}$
for all $\omega\in\Phi$. Moreover $\omega(\delta)=1$ if and only if
$\omega\in W\omega_\delta$. This implies $\delta\in V_0$. Let
$\ell_\delta:=\sum W\omega_\delta$ and
$\ell^\delta:=\ell-\ell_\delta$. These are $W$\_invariant linear
functions on $V$.

\Examples: 1. {\it Classical case:} Here
$\delta=(1,\ldots,1)$, $\omega_\delta(z)=z_n$,
$W\omega_\delta=\{z_1,\ldots,z_n\}$,
$\l_\delta(z)=\sum_iz_i=\ell(z)$, and $\ell^\delta(z)=0$.

\noindent
2. {\it Semiclassical case:} Here
$\delta=(1,\ldots,1)$, $\omega_\delta(z)=z_n$,
$W\omega_\delta=\{z_i\mid n-i\ {\rm even}\}$,
$\l_\delta(z)=\sum_{i:n-i\ \rm even}z_i$, and
$$\eqno{}
\ell^\delta(z)=\cases{0&if $n$ is odd,\cr
\sum_{i=1}^n(-1)^{i-1}z_i&if $n$ is even.\cr}
$$
\medskip\noindent
For $\lambda\in\Lambda_+$ let
$$\eqno{}
c_\lambda^{(\delta)}=
c_\lambda^{(\delta)}(\rho):=\prod\limits_{\omega\in W\omega_\delta}
(\omega(\rho)+k_\omega)_{\omega(\lambda)},$$
where $(a)_n:=a(a+1)(a+2)\ldots(a+n-1)$ is the Pochhammer symbol. Up
to a sign, this is just the contribution of $W\omega_\delta$ to
$f_\lambda(-\rho)$.  Now we renormalize $\pq_\lambda$ as follows:
$$\eqno{}
\qq^{(\delta)}_\lambda(z):=
{c_\lambda^{(\delta)}\over d_\lambda}\pq_\lambda(z)=
(-1)^{\ell(\lambda)}c_\lambda^{(\delta)}\qq_\lambda(z).
$$
Then the generalized {\it binomial formula\/} is:

\Theorem BinomialTheorem. Let $\delta\in\Sigma^\vee\cap V^W$. Then
$$\eqno{BinomialFormula}
\qq^{(\delta)}_\lambda(z+\delta)=\sum
\limits_{\mu\in\Lambda_+\atop\ell^\delta(\mu)=\ell^\delta(\lambda)}
p_\mu(\rho+\lambda)\ \qq^{(\delta)}_\mu(z)
\quad\hbox{for every }\lambda\in\Lambda_+.
$$

\Proof: To emphasize dependence on $\rho$ we will also write
$p_\lambda(z;\rho)$ etc. Let $\rho':=\rho+\1\delta$ with
$s\in\CC$. Then it follows from the definitions that
$$\eqno{}
p_\lambda(z;\rho')=p_\lambda(z-\1\delta;\rho)\qquad\hbox{and}\qquad
f_\lambda(z;\rho')=f_\lambda(z-\1\delta;\rho).
$$
Hence
$$\eqno{}
q_\lambda(z;\rho')=
{f_\lambda(\rho'+\lambda;\rho')\over
f_\lambda(-\rho';\rho')}p_\lambda(z;\rho')=
{f_\lambda(\rho+\lambda)\over
f_\lambda(-\rho-s\delta)}p_\lambda(z-\1\delta)=
{f_\lambda(-\rho)\over
f_\lambda(-\rho-s\delta)}q_\lambda(z-\1\delta)
$$
Since the contributions of $\omega\in\Phi\setminus W\omega_\delta$ and
$\alpha\in\Delta$ cancel out, we have
$$\eqno{}
{f_\lambda(-\rho)\over
f_\lambda(-\rho-s\delta)}=
{c_\lambda^{(\delta)}(\rho)\over c^{(\delta)}_\lambda(\rho+s\delta)}.
$$
Now we apply the transposition formula \cite{DualFormula} with $\rho'$
instead of $\rho$. We also replace $z$ by $z+\1\delta$. Then we
obtain:
$$\eqno{DualFormula2}
c^{(\delta)}_\lambda(\rho)q_\lambda(-z-s\delta)=
\sum_\mu(-1)^{\ell(\mu)}p_\mu(\rho+\lambda)
{c^{(\delta)}_\lambda(\rho+s\delta)\over
c^{(\delta)}_\mu(\rho+s\delta)}
\,c^{(\delta)}_\mu(\rho)q_\mu(z).
$$
Let $t$ be a formal parameter. In equation \cite{DualFormula2}, we
replace $z$, $s$ by $t^{-1}z$, $t^{-1}$, respectively, and
multiply by $t^{\ell(\lambda)}$. Thus, we get
$$\eqno{DualFormula3}
c^{(\delta)}_\lambda(\rho)t^{\ell(\lambda)}q_\lambda(-t^{-1}z-t^{-1}\delta)=
\sum_\mu p_\mu(\rho+\lambda)
A_\mu(t)\ [(-1)^{\ell(\mu)}c^{(\delta)}_\mu(\rho)t^{\ell(\mu)}q_\mu(t^{-1}z)].
$$
with
$$\eqno{}
A_\mu(t):=t^{\ell(\lambda-\mu)}
{c^{(\delta)}_\lambda(\rho+t^{-1}\delta)\over
c^{(\delta)}_\mu(\rho+t^{-1}\delta)}.
$$
Now, we take the limit for $t\pfeil0$. The left hand side of
\cite{DualFormula3} becomes $\qq^{(\delta)}_\lambda(z+\delta)$ while the
expression in brackets on the right hand side tends to
$\qq^{(\delta)}(z)$. Finally, we have
$$\eqno{}
(\omega(\rho+t^{-1}\delta)+k_\omega)_{\omega(\lambda)}=
(t^{-1}+\omega(\rho)+k_\omega)_{\omega(\lambda)}=t^{-\omega(\lambda)}+
\ldots
$$
where again ``$\ldots$'' means ``terms of higher order in $t$''.
Thus
$$\eqno{}
c^{(\delta)}_\lambda(\rho+t^{-1}\delta)=t^{-\ell_\delta(\lambda)}+\ldots
$$
and
$$\eqno{order}
A_\mu(t)=t^{\ell^\delta(\lambda-\mu)}+\ldots.
$$
By the Extra Vanishing \cite{ExtraVanishing} only those $\mu$ in
\cite{DualFormula3} have to be considered for which
$\tau:=\lambda-\mu\in\Lambda$. Thus the binomial formula
\cite{BinomialFormula} is proved when we show that
$\ell^\delta(\tau)\ge0$ for all $\tau\in\Lambda$.

Since $\ell^\delta$ is linear we may assume $\ell(\tau)=1$ since those
$\tau$'s generate $\Lambda$. Because $\ell^\delta$ is $W$\_invariant,
we may moreover assume that $\tau\in\Sigma^\vee$. Now consider formula
\cite{DualFormula3} with $\lambda=\tau$. Then the right-hand side has
only two non\_vanishing terms summands, corresponding to $\mu=\tau$
and $\mu=0$. Thus
$$\eqno{}
c^{(\delta)}_\tau(\rho)t q_\tau(-t^{-1}z-t^{-1}\delta)=
A_0(t)-A_\tau(t)c^{(\delta)}_\tau(\rho)t q_\tau(t^{-1}z)
$$
Since the limit $\|lim|_{t\pfeil0}A_0(t)$ exists and $q_\tau(z)$ is a
non\_constant polynomial of degree $1$ also
$\|lim|_{t\pfeil0}A_\tau(t)$ exists. Therefore
$\ell^\delta(\tau)\ge0$ by \cite{order}.\qed

\noindent
Putting $z=0$, we get as an immediate consequence an evaluation formula:

\Corollary. For all $\lambda\in\Lambda_+$ and $\delta\in\Sigma^\vee\cap V^W$
holds
$$\eqno{}
\qq^{(\delta)}_\lambda(\delta)=\cases{1&if
$\ell^\delta(\lambda)=0$,\cr0&otherwise.\cr}
$$

\Remark: Consider the classical case. Then the binomial formula
\cite{BinomialFormula} is due to Okounkov\_Olshanski, \cite{OlOk}.
Before that, Lassalle, \cite{Las}, used the binomial formula to {\it
define} the ``generalized binomial coefficients''
$p_\mu(\rho+\lambda)$. We see now that this was only possible because
$\ell^\delta=0$. For arbitrary multiplicity free actions,
Yan \cite{Yan} took another approach to define $p_\mu(\rho+\lambda)$
{}from the homogeneous polynomials $\pq_\lambda$, namely via the formula
$$\eqno{}
{1\over k!}\ell(z)^k\pq_\lambda(z)=
\sum_{\mu\in\Lambda_+\atop\ell(\mu-\lambda)=k}p_\lambda(\rho+\mu)\pq_\mu(z)
$$
which follows readily from \cite{Pieriell}. Yet another construction
can be found in \cite{BenRat}. Observe though, that none of these
approaches give the polynomiality nor the $W$\_invariance of
$p_\lambda$. Also the latter two constructions work only for those
$\rho\in V_0$ which actually come from a multiplicity free action.

\beginsection Beispiel. Example: the rank one case

In this section, we illustrate the main assertions of this paper with the
rank one case.

\medskip\noindent
{\it Section \cite{Fourier}:} Let $G=GL_n(\CC)$ and $U=\CC^n$, the
defining representation. Then $\cP^\lambda=S^\lambda(\CC^n)^\vee$, the
space of homogeneous polynomials of degree $\lambda\in\NN$. The
algebra of invariant differential operators is generated by
$\xi=\sum_ix_i{\partial\over\partial x_i}$, the Euler vector field. The
eigenvalue of $\xi$ on $\cP^\lambda$ is $\lambda$, hence
$c_\xi(z)=z$.

The parabolic $P$ is the stabilizer of the line $\CC
e_1\subseteq\CC^n$. Denote the weights of $\CC^n$ by
$\epsilon_i$. Then the roots in the unipotent radical of $P$ are
$\epsilon_1-\epsilon_i$. Thus,
$$\eqno{}
\rho={1\over2}\left(\sum_{i=2}^n(\epsilon_1-\epsilon_i)+
\sum_{i=1}^n\epsilon_i\right)={n\over2}\epsilon_1.
$$
Thus $p_\xi(z)=z-{n\over2}$. On the other hand we have
$$\eqno{}
\9\xi=\sum_{i=1}^n\left(-{\partial\over\partial x_i}\right)x_i
=-\sum_{i=1}^n(x_i{\partial\over\partial x_i}+1)=-\xi-n.
$$
Thus
$$\eqno{}
p_{\9\xi}(z)=p_{-\xi-n}(z)=-(z-{n\over2})-n=-z-{n\over2}=p_\xi(-z).
$$
\medskip
\noindent
{\it Section \cite{Capelli}:}
In the rank one case we have
$$\eqno{}
V=\CC,\quad W=1,\quad \Lambda_+=\NN,\quad\hbox{and}\quad\ell(z)=z.
$$
Moreover
$$\eqno{}
\Sigma^\vee=\Lambda_1=\{1\},\quad
\Phi=\Phi^+=\Sigma=\{z\},\quad\hbox{and}\quad\Delta=\emptyset.
$$
We have $V_0=\CC$ and put $\rho=s$. Thus every $\rho$ is non\_integral
while ``strongly dominant'' means $s\not\in-{1\over2}\NN$.

The polynomial $p_\lambda\in\cP=\CC[z]$ vanishes in
$z=s,s+1,\ldots,s+\lambda-1$ and is 1 in $z=s+\lambda$. There is
indeed only one such polynomial namely
$$\eqno{}
p_\lambda(z)={z-s\choose\lambda}.
$$
We have $f_\tau(z)=[z-s\downarrow\tau]$ for $\tau\in\NN$. Thus
$$\eqno{}
L=(z-s)T\quad\hbox{and}\quad E=z-(z-s)T=(z-s)\nabla+s
$$
where $T$ is the shift operator $T(h)(z)=h(z-1)$ and $\nabla:=1-T$.
Then an easy calculation shows
$$\eqno{E40}
D_h=\sum_{d=0}^\infty(-1)^d{z-s\choose d}(\nabla^dh)(z)\,T^d
\quad\hbox{for all}\quad h\in\CC[z].
$$
The equation $E(p_\lambda)=(s+\lambda)p_\lambda$ is equivalent to the
well known relation
$$\eqno{}
z\nabla{z\choose\lambda}=\lambda{z\choose\lambda}
$$
while $D_h(p_\lambda)=h(s+\lambda)p_\lambda$ gives, after using
\cite{E40} and some easy manipulations,
Newton's interpolation formula:
$$\eqno{E43}
h(x+z)=\sum_{d=0}^\infty{1\over d!}(\nabla^dh)(z)\,(x)_d
$$
(we substituted $s+\lambda=x+z$). Here $(x)_d=x(x+1)\ldots(x+d-1)$ is
the Pochhammer symbol. This can be used to rewrite formula
\cite{E40}. Since $T=1-\nabla$ we get
$$
\eqalignno{&D_h&=\sum_{0\le d\le m}(-1)^m{z-s\choose m}
(\nabla^mh)(z)\,(-1)^d{m\choose d}\nabla^d=\cr
E44&&=\sum_{d=0}^\infty{z-s\choose d}\left[\sum_{m=d}^\infty(-1)^{m-d}
{z-s-d\choose m-d}(\nabla^mh)(z)\right]\nabla^d\cr}
$$
If we apply $\nabla_z^m$ on both sides of \cite{E43} and then
substitute $x=s+d-z$ we get the expression in brackets of \cite{E44}. Thus
$$\eqno{E45}
D_h=\sum_{d=0}^\infty(\nabla^dh)(s+d){z-s\choose d}\nabla^d.
$$

\medskip\noindent
{\it Section \cite{Duality}:} According to \cite{E30} we have
$$\eqno{}
d_\lambda=(-1)^{\lambda}{[-2s\da\lambda]\over[\lambda\da\lambda]}
=(-1)^\lambda{-2s\choose\lambda}=
{2s-1+\lambda\choose\lambda}
$$
which affirms the evaluation formula \cite{E8}. Moreover, in the
geometric situation above with $GL_n(\CC)$ acting on $\CC^n$ we check
\cite{virdim}:
$$\eqno{}
\|dim|S^\lambda(\CC^n)^\vee={n-1+\lambda\choose\lambda}
$$
(since $s={n\over2}$). Furthermore,
$$\eqno{}
q_\lambda(z)={[z-s\da\lambda]\over[-2s\da\lambda]}
={(-z+s)_\lambda\over(2s)_\lambda}.
$$
Thus, the
transposition formula \cite{DualFormula} reads
$$\eqno{}
{(z+s)_\lambda\over(2s)_\lambda}=\sum_{\mu=0}^\lambda
(-1)^\mu{\lambda\choose\mu}{(-z+s)_\mu\over(2s)_\mu}.
$$
A direct proof boils down, after some manipulations, to the
Chu\_Vandermonde identity. Finally, the symmetry statement
\cite{Symmetrisch} becomes
$$\eqno{}
{(2s+\nu)_\lambda\over(2s)_\lambda}={(2s+\lambda)_\nu\over(2s)_\nu}
$$
which is easily verified directly.

\medskip\noindent
{\it Section \cite{Interpolation}:} The involutivity of the matrix
\cite{E41}
$$\eqno{}
\left((-1)^\mu{\lambda\choose\mu}\right)_{\lambda\mu}
$$
is well known. The transformation $h\mapsto\widehat h$ can be
rewritten as
$$\eqno{}
\widehat
h(s+\lambda)=\sum_{\mu=0}^\lambda(-1)^\mu{\lambda\choose\mu}h(s+\mu)=
(-1)^\lambda(\Delta^\lambda h)(s)
$$
where
$$\eqno{}
\Delta:=T^{-1}-1,\quad \hbox{i.e.},\quad(\Delta h)(z)=h(z+1)-h(z).
$$
Then the interpolation formula \cite{E11} becomes another
form of Newton interpolation (with $z=x+s$):
$$\eqno{}
h(x+s)=\sum_{\mu=0}^\infty(\Delta^\mu h)(s){x\choose\mu}
$$

\medskip\noindent
{\it Section \cite{Product}:} The scalar product \cite{E42} is
$$\eqno{}
\<g,h\>=\sum_{\mu=0}^\infty{2s-1+\mu\choose\mu}(\Delta^\mu
g)(s)(\Delta^\mu h)(s).
$$

\medskip\noindent
{\it Section \cite{sl-two}:} We have
$$
\eqalignno{&L^-&=-(z+s)T^{-1}\cr
&L^*&=(z-s)T-2z+(z+s)T^{-1}=z(\Delta-\nabla)+s(\Delta+\nabla)\cr}
$$
Since $\cB$ is the algebra generated by the $sl_2$\_triple
$$\eqno{}
\big((z-s)T,2z,-(z+s)T^{-1}\big)
$$
it is actually isomorphic to the universal enveloping algebra of
$sl_2(\CC)$.

\medskip\noindent {\it Section \cite{differential}:}\
We have $\pq_\lambda(z)={z^\lambda\over \lambda!}$ and
$\qq_\lambda(z)={(-1)^\lambda\over(2s)_\lambda}z^\lambda$. The algebra
$\cBq$ is generated by
$$\eqno{}
\Psi(L)=z,\ \Psi(E)=z{ d\over d z}+s,\ 
\Psi(L^*)=z{ d^2\over d z^2}+2s{ d\over d z}.
$$
Moreover, according to \cite{E45}:
$$\eqno{}
\Psi(D_h)=\Dq_h=\sum_{m=0}^\infty{1\over
m!}(\nabla^mh)(s+m)\,z^m{ d^m\over d z^m}.
$$

\medskip\noindent
{\it Section \cite{Binomial}:} We have $\delta=1$,
$c_\lambda^{(\delta)}=(2s)_\lambda$, and
$\qq_\lambda^{(\delta)}(z)=z^\lambda$. Thus, formula \cite{BinomialFormula} 
just specializes to the classical binomial formula
$$\eqno{}
(z+1)^\lambda=\sum_{\mu=0}^\lambda{\lambda\choose\mu}z^\mu.
$$

\hfuzz5pt

\beginrefs

\baselineskip11.5pt

\L|Abk:BenRat1|Sig:\\BR|Au:Benson, C., Ratcliff, G.|Tit:A
classification of multiplicity free actions%
|Zs:J. Algebra|Bd:181|S:152--186|J:1996||

\L|Abk:BenRat|Sig:\\BR|Au:Benson, C.; Ratcliff, G.%
|Tit:Combinatorics and spherical functions on the Heisenberg group%
|Zs:Represent. Theory {\rm(electronic)}|Bd:2|S:79--105|J:1998||

\L|Abk:Deb|Sig:De|Au:Debiard, A.|Tit:Polyn{\^o}mes de Tch{\'e}bychev et de
Jacobi dans un espace euclidien de dimension $p$|Zs:C.R. Acad. Sc. Paris%
|Bd:296|S:529--532|J:1983||

\L|Abk:HoUm|Sig:HU|Au:Howe, R., Umeda, T.|Tit:The Capelli
identity, the double commutant theorem, and multiplicity-free
actions|Zs:Math. Ann.|Bd:290|S:565--619|J:1991||

\L|Abk:Las|Sig:\\La|Au:Lassalle, M.|Tit:Une formule du bin\^ome
g\'en\'eralis\'ee pour les polyn\^omes de Jack|Zs:C. R. Acad. Sci.
Paris S\'er. I Math.|Bd:310|S:253--256|J:1990||

\L|Abk:La2|Sig:\\La|Au:Lassalle, M.|Tit:Coefficients binomiaux
g\'en\'eralis\'es et polyn\^omes de
Macdonald|Zs:J. Funct. Anal.|Bd:158|S:289--324|J:1998||

\L|Abk:Kac|Sig:Kac|Au:Kac, V.|Tit:Some remarks on nilpotent orbits%
|Zs:J. Algebra|Bd:64|S:190--213|J:1980||

\L|Abk:Annals|Sig:\\Kn|Au:Knop, F.|Tit:A Harish-Chandra homomorphism
for reductive group actions|Zs:Ann. of Math.
(2)|Bd:140|S:253--288|J:1994||

\Pr|Abk:Montreal|Sig:\\Kn|Au:Knop, F.|Artikel:Some remarks on
multiplicity free spaces|Titel:Proc. NATO Adv. Study Inst. on
Representation Theory and Algebraic Geometry|Hgr:A.~Broer,
G.~Sabi\-dussi, eds.|Reihe:Nato ASI Series C|Bd:514%
|Verlag:Kluwer|Ort:Dortrecht|S:301--317|J:1998||

\L|Abk:Semi|Sig:\\Kn|Au:Knop, F.|Tit:Semisymmetric polynomials
and the invariant theory of matrix vector pairs|Zs:{\tt math.RT/9910060}%
|Bd:-|S:26 pages|J:1999||

\L|Abk:Existence|Sig:\\Kn|Au:Knop, F.|Tit:Construction of commuting
difference operators for multiplicity free spaces|Zs:{\tt
math.RT/0006004}|Bd:-|S:28 pages|J:2000||

\L|Abk:SymCap|Sig:KnSa|Au:Knop, F.; Sahi, S.|Tit:Difference
equations and symmetric polynomials defined by their
zeros|Zs:Internat. Math. Res. Notices|Bd:10|S:473--486|J:1996||

\L|Abk:Leahy|Sig:Le|Au:Leahy, A.|Tit:A classification of
multiplicity free representations|Zs:J. Lie Theory|Bd:8|S:367--391|J:1998||

\L|Abk:Ok|Sig:Ok|Au:Okounkov, A.|Tit:Binomial formula for Macdonald
polynomials and applications|Zs:Math. Res. Lett.|Bd:4|S:533--553|J:1997||

\L|Abk:OO1|Sig:\\OO|Au:Olshanski, G.; Okounkov, A.|Tit:Shifted Schur
functions|Zs:St.~Petersburg Math. J.|Bd:9|S:73--146|J:1998||

\L|Abk:OlOk|Sig:\\OO|Au:Olshanski, G.; Okounkov, A.|Tit:Shifted Jack
polynomials, binomial formula, and applications|Zs:Math. Res. Lett.%
|Bd:4|S:69--78|J:1997||

\Pr|Abk:Sahi|Sig:Sa|Au:Sahi, S.|Artikel:The spectrum of certain
differential operators associated to symmetric space%
|Titel:Lie Groups and Geometry|Hgr:J.-L. Brylinski et al.,
eds.|Reihe:Progr. Math.%
|Bd:123|Verlag:Boston|Ort:Birkh{\"a}user|S:569--576|J:||

\L|Abk:Sek|Sig:Se|Au:Sekiguchi, J.|Tit:Zonal spherical functions on
some symmetric spaces|Zs:Publ. RIMS, Kyoto Univ.|Bd:12|S:455--459|J:1977||

\L|Abk:Up|Sig:Up|Au:Upmeier, H|Tit:Toeplitz operators on bounded symmetric
domains|Zs:Trans. Amer. Math. Soc.|Bd:280|S:221--237|J:1983||

\L|Abk:Yan|Sig:Yan|Au:Yan, Zhi Min|Tit:Special functions associated
with multiplicity free representations|Zs:Preprint|Bd:-|S:-|J:1992||

\endrefs

\bye